\documentclass[3p,review]{elsarticle}
\usepackage{amsmath,amssymb}
\usepackage{verbatimbox}
\usepackage{graphicx}
\graphicspath{{./Pics/}}
\usepackage{listings}
\usepackage[colorlinks=true,linkcolor=blue,citecolor=blue,urlcolor=blue]{hyperref}

\usepackage{siunitx}

\usepackage{enumitem}

\newcounter{hintcounter}

\newcommand{\hint}[1]{

  \begin{center}
    \emph{``Hint \thehintcounter:  {#1}''}
  \end{center}
  \stepcounter{hintcounter}
}

\newcommand{\N}{\mathcal{N}}
\newcommand{\T}{\mathcal{T}}
\newcommand{\D}{\mathcal{D}}
\newcommand{\Cl}{\mathcal{C}}
\newcommand{\I}{\mathrm{i}}
\newcommand{\Es}{\mathcal{S}}
\newcommand{\sphere}{{\mathbb{S}^2}}

\newcommand{\bx}{{\bf x}}
\newcommand{\by}{{\bf y}}
\newcommand{\bn}{{\bf n}}
\newcommand{\bz}{{\bf z}}
\newcommand{\bs}{{\bf s}}
\newcommand{\bd}{{\bf d}}
\newcommand{\bT}{{\bf T}}
\newcommand{\bN}{{\bf N}}
\newcommand{\bS}{{\bf S}}
\newcommand{\bD}{{\bf D}}
\newcommand{\de}{\mathrm{d}}

\newcommand{\meshtohrtf}{Mesh2HRTF{}}
\newcommand{\numcalc}{NumCalc{}}

\begin{document} 
\begin{frontmatter}{}
\title{NumCalc: An open-source BEM code for solving acoustic scattering problems}
\author[ARI]{W. Kreuzer\corref{CORRAUTH}}
\ead{wolfgang.kreuzer@oeaw.ac.at}
\author[ARI]{K. Pollack}
\ead{katharina.pollack@oeaw.ac.at}
\author[TUB]{F. Brinkmann}
\ead{fabian.brinkmann@tu-berlin.de}
\author[ARI]{P. Majdak}
\ead{piotr.majdak@oeaw.ac.at}

\affiliation[ARI]{
  organization = {Austrian Academy of Sciences, Acoustics Research Institute},
  addressline = {Dominikanerbastei 16, 3. Stock},
  postcode = {1010},
  city = {Vienna},
  country = {Austria}
}

\affiliation[TUB]{
  organization = {Audio Communication Group, Technical University of Berlin},
  addressline = {Straße des 17. Juni 135},
  postcode = {10623},
  city = {Berlin},
  country = {Germany}
}

\cortext[CORRAUTH]{Corresponding author.}

\begin{abstract}
The calculation of the acoustic field in or around objects is an important task in acoustic engineering. The open-source project \meshtohrtf{} and its BEM core \numcalc{} provide users with a collection of free tools for acoustic simulations without the need of having an in-depth knowledge into numerical methods. However, we feel that users should have a basic understanding with respect to the methods behind the software they are using. We are convinced that this basic understanding helps in avoiding common mistakes and also helps to understand the requirements to use the software. To provide this background is the first motivation for this article. A second motivation for this article is to demonstrate the accuracy of \numcalc{} when solving benchmark problems. Thus, users can obtain an idea about the accuracy as well as requirements on the memory and CPU requirements when using NumCalc. Finally, this article provides detailed information about some aspects of the actual implementation of BEM that are usually not mentioned in literature, e.g., the specific version of the fast multipole method and its clustering process or how to use frequency-dependent admittance boundary conditions.
\end{abstract}

\begin{keyword}
  BEM \sep Software \sep Fast Multipole Method
\end{keyword}
\journal{Engineering Analysis with Boundary Elements}
\end{frontmatter}{}
\section{Introduction}\label{Sec:Introduction}
The boundary element method (BEM, ~\cite{chen2008,SauSch10}) has a long tradition for numerically solving the Helmholtz equation in 3D that models acoustic wave propagation in the frequency domain. Compared to other methods such as the finite element method~\cite{MarNol08} (FEM), the BEM has the advantage that only the surface of the scattering object needs to be considered and that for external problems the decay of the acoustic wave towards infinity is already included in the formulation. This renders the BEM an attractive tool for calculating sound-wave propagation. However, the BEM has the drawback that the generated linear system of equations has a densely populated system matrix. To tackle this problem,  the combination of BEM and fast methods for matrix-vector multiplications such as the fast multipole method (FMM,~\cite{Coifmanetal93,Rokhlin90,Rokhlin93}) or $\mathcal{H}$-matrices~\cite{Hackbusch15} were introduced. Eventually, BEM-based acoustic modeling became an increasingly popular tool for engineers, as illustrated by various software packages such as COMSOL~\cite{comsol}, FastBEMAcoustics~\cite{fastbem}, AcouSTO~\cite{Acousto}, BEMpp~\cite{Betcke2021}, or OpenBEM~\cite{Henrquez2010OpenBEMA}. 

\numcalc{}~\cite{Kreuzeretal22a} is an open-source program written in C$++$ for solving the Helmholtz equation in 3D. It is based on collocation with constant elements in combination with an implementation of the FMM to speed up matrix-vector multiplications. \numcalc{} is distributed in the framework of 
\meshtohrtf{}~\cite{Fabianspaper,mesh2hrtf,Ziegelwangeretal15}, which is aimed at the numerical calculation and post processing of head-related transfer functions (HRTFs,~\cite{moller_head-related_1995}). HRTFs describe the direction-dependent filtering of a sound source due to the listener's body especially by the head and ears and are usually acoustically measured~\cite{majdak2007,pollack2022modern}.

The basic idea for \meshtohrtf{} is to provide an open-source code that is easy-to-use for researchers without the need for an extensive background in mathematics or physics, specifically targeting users not familiar with BEM. To this end, \meshtohrtf{} contains an add-on for the open-source program Blender \cite{blendermanual} and automatically generates the input for \numcalc. The accompanying project websites \cite{mesh2hrtf,mesh2hrtfwiki} provide tutorials for creating a valid description of the geometry of the scattering object and its materials. Despite the encapsulation of the BEM calculations from the typical \meshtohrtf{} user, some rules are required for correct calculations and their compliance cannot be ensured by the \meshtohrtf{} interface. Hence, the main goal of the article is to explain such rules, provide insights to what can happen if they are violated, and how to avoid common mistakes. Throughout the article we point out the most important rules as hints:
\hint{Look out for these lines.} 

Although part of the \meshtohrtf{}-project, \numcalc{} is \emph{not restricted} to calculations of HRTFs only. The input to \numcalc{} is a general definition of surface meshes based on triangular or plane quadrilateral elements. Arbitrary admittance boundary conditions can be assigned to each element and external sound sources can be modelled as plane waves or as point sources (see \ref{Sec:Mesh} and \ref{App:Input}). \numcalc{} can be used to solve general wave scattering problems in 3D using the BEM.

This article describes the general features of \numcalc{} in detail for the first time in one manuscript. It provides all necessary information for users to responsibly use \numcalc{} for their problem at hand, including some benchmark tests to measure its accuracy and needs with respect to memory consumption. Besides  specific aspects of the implementation of the BEM, such as the quasi-singular quadrature or the possibility to define frequency-dependent boundary conditions, we for the first time give a detailed description on the specific flavour and implementation of the FMM. This includes a description, how the clusters are generated and how one can derive all matrices necessary for the fast multipole method. Compared to other standard implementations of the FMM the root level already consists of multiple clusters and the local element-to-cluster expansion matrices on \emph{all} levels of a multilevel FMM are stored explicitly, which reduces computing time. As the root level already consists of many clusters, the truncation length of the multipole expansion and in turn the memory necessary to use NumCalc can be kept relatively low. Finally, we provide numerical experiments that  illustrate the accuracy of \numcalc{}, show how the mesh is clustered in the default setting and discusses the memory required for the FMM algorithm based on the clustering. 

This article is structured as follows: Sec.~\ref{Sec:BriefBEM} describes general aspects of the BEM and the way they are implemented in \numcalc. Sec.~\ref{Sec:Implementation} provides specific details of the implementation of the quadrature and clustering in \numcalc{}. Sec.~\ref{Sec:Critical} focuses on the critical components of the implementation and their effect on accuracy and computational effort. Sec.~\ref{Sec:Benchmark} presents three benchmark problems solved with \numcalc{} to demonstrate the achievable accuracy and the need of computing power. Finally, we summarize and discuss the results in Sec.~\ref{Sec:Discussion}. To keep the article self contained, \ref{Sec:Mesh} provides the definition of a mesh to describe the geometry of the scatterer. In \ref{App:Memory} a method to estimate the computer memory needed by \numcalc{} is presented, \ref{App:Input} provides examples for \numcalc{} input files including how to define frequency-dependent admittance boundary conditions, \ref{App:Commandline} gives a description of possible command line parameters, and \ref{App:Output} lists the different output files generated by \numcalc.
\section{BEM and Its Implementation in \numcalc}\label{Sec:BriefBEM}
In a nutshell, \numcalc{} implements the direct BEM using collocation with constant elements in combination with the Burton-Miller method. By default, an iterative solver is used to solve the system of equations and the matrix-vector multiplications needed by this solver can be made efficiently with the FMM. To keep the article self-contained, we start with a brief introduction to collocation BEM and explain the requirements for the meshes used for the BEM. 

\begin{figure}
  \begin{center}
  \includegraphics[width=0.3\textwidth]{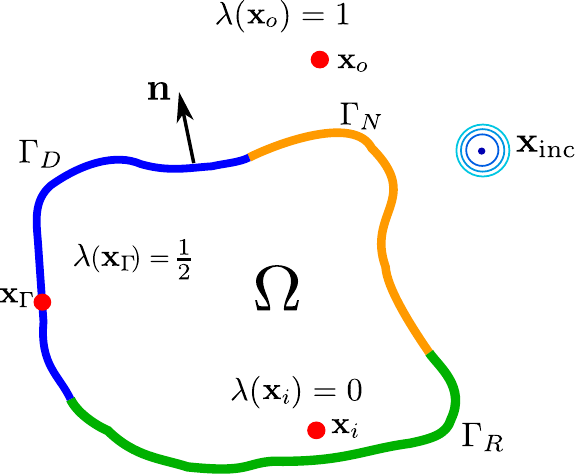}
  \caption{Scheme of an exterior problem with scatterer $\Omega$ with boundary $\Gamma = \Gamma_D \cup \Gamma_N \cup \Gamma_R$, a point-source at $\bx_{\text{inc}}$, interior point $\bx_i$, exterior point $\bx_o$ and a node on the boundary $\bx_{\Gamma}$.}\label{Fig:Schemedrawing}
  \end{center}
\end{figure}

Figure~\ref{Fig:Schemedrawing} shows a typical situation of an (exterior) acoustics wave scattering problem. An incoming point source outside the scatterer $\Omega$ is positioned at a point $\bx_{\mathrm{inc}}$. The acoustic wave produced by this source is scattered and reflected at the surface $\Gamma$ of the scatterer. With the BEM, the acoustic field on the surface $\Gamma$ at nodes $\bx_\Gamma$ as well as on exterior nodes $\bx_{\text{o}}$ outside $\Gamma$ is determined. The main goal of \numcalc{} is to derive the acoustic pressure $p$ and the particle velocity $v$ at both the points on the surface $\Gamma$ of the scatterer $\Omega$ and the evaluation points $\bx_{\text{o}}$ in the free-field $\mathbb{R}^3\setminus \Gamma$. To this end, \numcalc{} solves the Helmholtz equation describing the free-field sound propagation in the frequency domain:

\begin{equation}\label{Equ:Helmholtz}
\nabla^2 p(\bx,\omega) + k^2 p(\bx,\omega) = p_0(\bx,\bx_{\text{inc}},\omega),\quad \bx \in \mathbb{R}^3\setminus \Gamma,
\end{equation}
with $p$ being the sound pressure, $k = \frac{\omega}{c} = \frac{2\pi f}{c}$ being the wavenumber for a given frequency $f$ and speed of sound $c$, and $p_0$ being an external sound source located at $\bx_{\text{inc}}$. 

To model the acoustic properties of the surface $\Gamma$, a boundary condition (bc) is defined for each point on the surface $\Gamma$. Either the sound pressure $p$ (Dirichlet bc), the particle velocity $v$ (Neumann bc) or a combination of both (Robin bc) need to be known:
\begin{align*}
  p(\bx) &= f(\bx)\qquad  \bx \in \Gamma_D,\\
  v(\bx) & = g(\bx)\qquad \bx \in \Gamma_N,\\
  \alpha p(\bx) + \beta v(\bx) &= h(\bx)\qquad \bx \in \Gamma_R,
\end{align*}
for given functions $f,g,$ $h$, and constants $\alpha,\beta$. 
\subsection{The boundary integral equation (BIE)}
The acoustic field at a point $\bx$ can be represented by the sound pressure $p(\bx)$ and the particle velocity $v(\bx)$ at this point. Both can be derived in the frequency domain from the velocity potential $\phi(\bx)$ using
\begin{equation}\label{Equ:utopandv}
  p(\bx) = \I \omega \rho \phi(\bx),\quad v(\bx) = \frac{\partial \phi}{\partial n_x}\left(\bx\right) := \nabla \phi(\bx) \cdot \bn(\bx),
\end{equation}
where $\omega = 2\pi f$ is the angular frequency, and $\rho$ denotes the density of the medium. If $\bx \in \Gamma$, the particle velocity is defined using the normal vector $\bn(\bx)$ to the surface at $\bx$ pointing to the outside. If the point $\bx$ lies in the free field, it can be any arbitrary vector with $||\bn|| = 1$.

In order to solve the Helmholtz equation, the partial differential equation needs to be transformed into an integral equation. In general, there are two possibilities for this transformation, direct and indirect BEM (cf.~\cite{SauSch10}). In \numcalc{}, the direct formulation based on Green's representation formula is used:
\begin{equation}\label{Equ:BIE_I}
  \lambda(\bx) \phi(\bx) - \tau \int\limits_\Gamma H(\bx,\by) \phi(\by) \de\by + \tau \int\limits_\Gamma G(\bx,\by) v(\by) \de\by = \phi_{\text{inc}}(\bx),
\end{equation}
where $G(\bx,\by) = \frac{e^{\I k ||\bx - \by||}}{4\pi ||\bx - \by||}$ is the Green's function for the free-field Helmholtz equation, $H(\bx,\by) = \frac{\partial G}{\partial n_y}(\bx,\by) = \nabla G(\bx,\by) \cdot \bn(\by) = G(\bx,\by)\left(\I k - \frac{1}{r}\right) \frac{\by - \bx}{r}\cdot \bn(\by)$, $r = ||\bx - \by||$, and $\phi_{\text{inc}}$ is the velocity potential caused by an external sound source, e.g., a point source positioned at $\bx_{\mathrm{inc}}$. The parameter $\tau$ indicates the problem setting. If $\tau = -1$, the acoustic field \emph{inside} the object is of interest, in this case, we speak of an ``interior'' problem. If $\tau = 1$, the acoustic field outside the object is of interest, in that case, we have an ``exterior'' problem. The parameter $\lambda(\bx)$  in Eq.~(\ref{Equ:BIE_I}) depends on the position of the point $\bx$ (see also Fig.~\ref{Fig:Schemedrawing}). If $\bx $ lies in the interior of the object $\Omega$, $\lambda(\bx) = 0$. If $\bx$ lies outside the object, $\lambda(\bx) = 1$. If $\bx$ lies on a smooth part of the surface, $\lambda(\bx) = \frac12$.

The integrands in Eq.~(\ref{Equ:BIE_I}) depend on the derivative with respect to the normal vector $\bn_y := {\bf n}(\by)$ pointing away from the scatterer. Normal vectors pointing inside the object cause a sign change and the roles of interior and exterior problems are switched. As a consequence, Eq.~(\ref{Equ:BIE_I}) may produce incorrect results. To avoid this, it is essential that $\bn_y$ points away from the scatterer, however, \numcalc{} does not check this. This check is left to the user: \hint{Check if your normal vectors point away from the scatterer.}

Note that if the distance $r = ||\bx - \by||$ between two points $\bx$ and $\by$ becomes small, $G$ and $H$ behave like $\frac{1}{r}$. Thus, special methods are needed to numerically calculate the integrals over $\Gamma$. We describe these methods in Sec.~\ref{Sec:Quadrature}. 

Note also that Eq.~(\ref{Equ:BIE_I}) cannot be solved uniquely at ``[...] `resonant' wavenumbers (or eigenvalues) for a related interior problem''~\cite{BurMil71}. For exterior problems, this non-uniqueness problem can be solved using either CHIEF points~\cite{Schenck68} or the Burton-Miller method~\cite{BurMil71}. \numcalc{} implements the Burton-Miller method and the implementation is described in Sec.~\ref{Sec:Burton-Miller}.
\subsection{Discretization of the Geometry and Shapefunctions}
To transform the BIE into a linear system of equations, the unknown solution of Eq.~(\ref{Equ:BIE_I}) on the surface (either $\phi(\bx)$ or $v(\bx)$ depending on the boundary condition at \bx)  needs to be approximated by simple ansatz functions (= shape functions). For that matter, the geometry of the scatterer's surface is approximated by a mesh consisting of $N$ elements $\Gamma_j, j = 1,\dots,N$.  In \numcalc{}, these elements can be triangular or \emph{plane} quadrilateral, and on each of them the BIE solution is assumed to be constant, e.g.:
$$
\phi(\bx) = \left\{
  \begin{array}{cc}
    \phi_j & \bx \in \Gamma_j,\\
    0 & \text{otherwise}.
  \end{array}
\right.
$$
The number of elements $N$ is a compromise between calculation time and numerical accuracy. The numerical accuracy is actually frequency dependent because the unknown solution with its oscillating parts needs to be represented accurately by the constant ansatz functions. As a rule of thumb, the mesh should contain six to eight elements per wavelength~\cite{marburgsix2002}. \numcalc{} does a check on that and displays a warning if the average size of the elements is too low, however, it will still do the calculations. This is because, in certain cases, acceptable results can be achieved, when ``irrelevant'' parts of the mesh have a coarser discretization~\cite{Ziegelwangeretal16,Palm2021a}. The decision on the number of elements is within the responsibility of the user.
\hint{Check if your mesh contains 6-to-8 elements per wavelength in the relevant geometry regions.}
It is also known that for some problems and especially close to sharp edges and corners a finer discretization than 6-to-8 elements per wavelength may be necessary~\cite{Kreuzer22}. One should be aware that the 6-to-8 elements per wavelength rule is just a rule of thumb (see the discussion about this rule and the resulting error in \cite{marburgsix2002}), that works well in practice. If users are unsure about the accuracy of their solution, it is suggested to use some a-posteriori error estimator, cf.~\cite{Kurzetal21}, or to repeat the calculation with a mesh with different element length and to compare the result with the solution of the original mesh (see also the experiments in Section~\ref{Sec:Errororder}). In general, this approach may not be a feasible for a lot of practical applications that already involve meshes with a large number of elements.

The elements of the mesh are also used to calculate the integrals over the surface $\Gamma$. In \numcalc{}, integrals over $\Gamma$ are replaced by sums of integrals over each element, for example: 
\begin{equation}\label{Equ:IntI}
\int\limits_{\Gamma} H(\bx,\by)\phi(\by) \de\by \approx  \sum_{j=1}^N \int_{\Gamma_j} H(\bx,\by) \phi(\by) \de\by \approx \sum_{j=1}^N  \left( \int\limits_{\Gamma_j} H(\bx,\by) \de\by \right)\phi_j.
\end{equation}
As $\Gamma_j$ represents a simple geometric element, the integrals over each element can be calculated relatively easily. By using the shape functions in combination with the boundary conditions, the unknown solution of Eq.~(\ref{Equ:BIE_I}) can be reduced to $N$ unknowns (either $\phi_j$ or $v_j$ depending on the given boundary condition).
\subsection{Collocation}\label{Sec:Collocation}
In order to derive a linear system of equation, \numcalc{} uses the collocation method, i.e., after the unknown solution is represented using the shape functions, Eq.~(\ref{Equ:BIE_I}) is evaluated at $N$ collocation nodes $\bx_i$ for $i = 1,\dots,N$. For constant elements these nodes are, in general, the midpoints of each element $\Gamma_i$ defined as the mean over the coordinates of the element's vertices. Note that because the collocation nodes are located inside the smooth plane element, $\lambda(\bx_i) = \frac12$ for $i = 1,\dots, N$, and Eq.~(\ref{Equ:BIE_I}) becomes
\begin{equation}\label{Equ:BIESystem}
\frac{\phi_i}{2} - \tau \sum_{j=1}^N \left(\int_{\Gamma_j} H(\bx_i,\by)\de\by\right) \phi_j
+ \tau \sum_{j=1}^N \left(\int_{\Gamma_j} G(\bx_i,\by) \de\by \right)v_j = \phi_{\mathrm{inc}} (\bx_i),
\end{equation}
where $\phi_i = \phi(\bx_i)$ and $v_i = v(\bx_i),\; i = 1,\dots, N,$ are the velocity potential and the particle velocity at the midpoint $\bx_i$ of each element, respectively. By defining the vectors
$$
{\pmb \phi} := (\phi_1,\dots,\phi_N)^\top,\quad {\bf v} := ({v_1,\dots,v_N})^\top,\quad \pmb{\phi}_{\mathrm{inc}}=(\phi_{\mathrm{inc}}(\bx_1),\dots,\phi_{\mathrm{inc}}(\bx_N))^\top,
$$
and the matrices
\begin{equation}\label{Equ:Matrices_I}
{\bf H}_{ij} = \int\limits_{\Gamma_j} H(\bx_i,\by) \de\by, \quad {\bf G}_{ij} = \int\limits_{\Gamma_j}G(\bx_i,\by)\de\by
\end{equation}
we obtain the discretized BIE in matrix-vector form:
\begin{equation}\label{Equ:LinSystem}
\frac{\pmb\phi}{2} - \tau {\bf H\pmb\phi} + \tau {\bf G v} = {\pmb\phi}_{\mathrm{inc}},
\end{equation}
which is a basis for further consideration to solve the Helmholtz equation numerically.
\subsection{Incoming Sound Sources}\label{Sec:Incoming}
Incoming sound sources $\phi_{\mathrm{inc}}$ are the sources not positioned on the surface of the scatterer but located in the free field. \numcalc{} implements two types of incoming sound sources: Plane waves and point sources. A plane wave is defined by its strength $S_0$ and its direction $\bd\in \mathbb{R}^3, ||\bd|| = 1$:
$$
\phi_{\mathrm{inc}}(\bx) = S_0 e^{\I k \bx \cdot \bd}, 
$$
and a point source is defined by its strength $S_0$ and its position $\bx^*$:
\begin{equation}\label{Equ:PointSrc}
\phi_{\mathrm{inc}}(\bx) = S_0 \frac{e^{\I k ||\bx - \bx^*||}}{4\pi||\bx - \bx^*||},
\end{equation}
with $k$ being the wavenumber. 
From the definition of the point source it becomes clear that if the source position $\bx^*$ is very close to the surface, the denominator in Eq.~(\ref{Equ:PointSrc}) will become zero causing numerical problems. Thus, \hint{Have at least the distance of one average edge length between an external sound source and the surface.}
If sound sources need to be positioned close to the surface or even on the surface, a velocity boundary condition different from zero at the elements close to the source may be a better choice.
\subsection{The Burton-Miller Method}\label{Sec:Burton-Miller}
\numcalc{} implements the Burton-Miller approach to ensure a stable and unique solution of Eq.~(\ref{Equ:LinSystem}) for exterior problems at all frequencies. The final system of equations is derived by combining Eq.~(\ref{Equ:BIE_I}) and its derivative with respect to the normal vector at $\bx$
\begin{equation}\label{Equ:dBIE}
\frac{v(\bx)}{2} - \tau \int_{\Gamma} E(\bx,\by) \phi(\by) \de\by + \tau \int_{\Gamma} H'(\bx,\by) v(\by) \de\by = v_{\mathrm{inc}}(\bx), 
\end{equation}
where $E(\bx,\by) = \frac{\partial^2 G(\bx,\by)}{\partial n_x \partial n_y}$ and $H'(\bx, \by) = \frac{\partial G(\bx,\by)}{\partial n_x}$. The matrix-vector representation of Eq.~(\ref{Equ:dBIE}) is then derived in a similar manner as for Eq.~(\ref{Equ:LinSystem}). 

The final system of equations to be solved reads as
\begin{equation}\label{Equ:SystemI}
  \frac{\left( \pmb\phi - \gamma \bf v\right)}{2} - \tau  \left( {\bf H - \gamma E}\right)\pmb{\phi} + \tau \left( {\bf G - \gamma H'}\right){\bf v} = \pmb{\phi}_{\text{inc}} - \gamma {\bf v}_{\text{inc}},
\end{equation}
where 
\begin{equation}\label{Equ:Theintegrals}
  {\bf E}_{ij} =   \int\limits_{\Gamma_j}E(\bx_i,\by) \de\by,\; {\bf H}'_{ij} = \int\limits_{\Gamma_j}H'(\bx_i,\by) \de\by,\; {\bf v}_{\mathrm{inc}} = (v_{\mathrm{inc}}(\bx_1),\dots, v_{\mathrm{inc}}(\bx_N))^\top,
\end{equation}
and $\gamma := \frac{i}{k}$ is the coupling factor introduced by the Burton-Miller method \cite{BurMil71,Marburg14}. 

Because a boundary condition (pressure, velocity, or admittance) is given for each individual element, either $\phi_i$, $v_i$, or a relation between them is known and the final system of equations Eq.~(\ref{Equ:SystemI}) has only $N$ unknowns. As there are as many collocation nodes as unknowns, the linear system can be solved by regular methods, e.g., either by a direct solver~\cite{Andersonetal99} or by an iterative solver~\cite{Saad03}. For most problems, the matrices in that system are densely populated and the number of unknowns can be large. Iterative solvers are thus the most common option. However, for large $N$, the cost of the matrix-vector multiplications needed by an iterative solver becomes too high. For those systems, the BEM can be coupled with methods for fast matrix-vector multiplications, e.g., $\mathcal{H}$-matrices~\cite{Hackbusch15} or the FMM~\cite{Coifmanetal93}. 
\section{Specific Details on Deriving and Solving the System Equation}\label{Sec:Implementation}
\subsection{Quadrature}\label{Sec:Quadrature}
In order to solve the BIE, the integrals in ${\bf G, H, H'}$, and ${\bf E}$ defined in Eqs.~(\ref{Equ:Matrices_I}) and (\ref{Equ:Theintegrals}) need to be calculated numerically. If the collocation node $\bx_i$ lies in the element for which the integral needs to be calculated, the integrand becomes singular because of the singularities of the Green's function and its derivatives. In this case, these integrals have to be calculated using special methods for the singular and the hypersingular integrals. For integrals over elements that do not contain the collocation node, we distinguish between elements that are close to the collocation node and elements far away which are solved by the quasi-singular quadrature and the regular quadrature, respectively. 
\subsubsection{Singular Quadrature}\label{Sec:SingQuad}
For $i=j$, the integrands in Eq.~(\ref{Equ:Matrices_I}) behave like $\frac{1}{r}$ for $r:=||\bx_i - \by|| \rightarrow 0$. In this case, \numcalc{} regularizes the integrals by using a quadrature method based on~\cite{Duffy82}. If the element $\Gamma_i$ is triangular, it is subdivided into six triangles, where the collocation node $\bx_i$, at which the singularity occurs, is a vertex of each of the sub-triangles. For a quadrilateral element, eight triangular sub-elements are created. Similar to~\cite{Duffy82}, the integrals over each of these sub-triangles are transformed into integrals over the unit square which removes the singularity. After that transformation, the integral is computed using a Gauss quadrature with $4\times 4$ quadrature nodes. 
For the singular quadrature with kernel ${\bf H}(\bx,\by)$, we can use the fact that, if $\bx$ and $\by$ are in the same planar element, $(\bx - \by)$ is orthogonal to $\bn_y$. This yields $H(\bx,\by) = \nabla_y G(\bx,\by)\cdot \bn_y = \frac{d G}{dr}\frac{\partial r}{\partial \by} \cdot \bn_y = -\frac{d G}{dr} \frac{\bx - \by}{r} \cdot \bn_y = 0$. The same argument can be used for the contributions with respect to ${\bf H}'$.
\subsubsection{Hypersingular Quadrature}
To calculate the integrals involving the hypersingular kernel $E(\bx,\by)$, the integral involving the hypersingular kernel is split into integrals over the edges of the element and the element itself~\cite{Krishnasamyetal90}:
\begin{equation}
{\bf E}_{ii} = \int\limits_{\Gamma_i}E(\bx_i,\by)\de\by= \sum_{\ell,m,n = 1}^3 \bn_\ell \epsilon_{\ell m n} \oint\limits_{\partial {\Gamma_i}}\frac{\partial G}{\partial y_m}(\bx_i,\by) dy_n + k^2\int\limits_{\Gamma_i}G(\bx_i,\by) \de\by, %= \frac{\partial G}{\partial n_y} \otimes 
\end{equation}
where $\epsilon_{\ell mn}$ is the Levi-Civita-symbol
$$
\epsilon_{\ell mn} = \left\{
  \begin{array}{rl}
    1 & \text{ if } (\ell,m,n) \in \{ (1,2,3), (2,3,1),(3,1,2) \},\\
    -1 & \text{ if } (\ell,m,n) \in \{ (3,2,1), (1,3,2), (2,1,3) \},\\
    0 & \text{ if } \ell=m, \ell=n, m=n,
  \end{array}
\right.
$$
and $\by = (y_1,y_2,y_3)$ is a point on the element $\Gamma_i$, $\bn_\ell$ is the $\ell$-th component of the normal vector $\bn$, and $\partial {\Gamma_i}$ denotes the edges of the element.
For the integral over $\partial \Gamma_i$ each edge is subdivided into four equally sized parts and on each part a Gaussian quadrature with 3 nodes is used. The integral over the element $\Gamma_i$ is calculated using the singular quadrature described in Sec.~\ref{Sec:SingQuad}.
\subsubsection{Quasi-singular Quadrature}
\begin{figure}[!h]
  \begin{center}
  \includegraphics[width=0.6\textwidth]{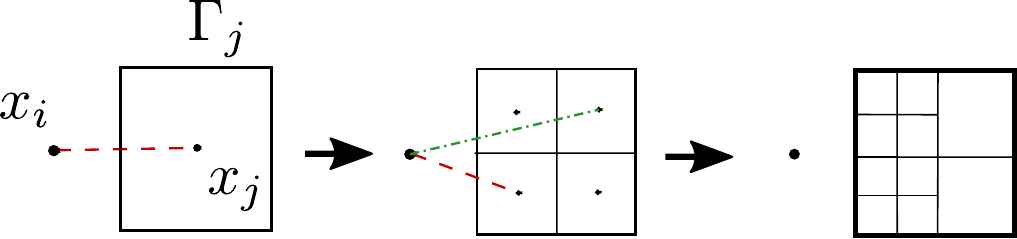}
  \caption{Example of two levels of element subdivisions in the quasi-singular quadrature and a quadrilateral element. The dashed (red) lines indicate that the distance between the two nodes is too small compared to the average edge length. In this case, the element needs to be divided. The dashed-dotted (green) line indicates that the distance is sufficiently large such that no further division is necessary.}\label{Fig:Subdivide}
  \end{center}
\end{figure}

\numcalc{} calculates the quasi-singular integrals in two steps: First, if for a given collocation node $\bx_i$ and an element $\Gamma_j$ with area $A_j$ and midpoint $\bx_j$
$$
\tilde{r} = \frac{||\bx_i - \bx_j||}{\sqrt{A_j}} < 1.3,
$$
$\Gamma_j$ is subdivided into four subelements. This procedure is repeated for all subelements until the condition is fulfilled. As a result, the distance between the collocation node and the subelement will be larger than the edge length of the subelement. With this procedure the subelements become smaller close to the collocation node (see, e.g., Fig.~\ref{Fig:Subdivide}). The factor 1.3 was heuristically chosen as a good compromise between accuracy and efficiency. Once the subelements are constructed, the integrals over them can be solved using the regular quadrature.

Note that the distance between the element midpoints is just an approximation for the actually required smallest distance between the collocation node and \emph{all} the points of an element. \numcalc{} uses the distance between the element midpoints because calculations can be done efficiently and it works well for regular elements. Thus,
\hint{Keep the elements as regular as possible.}
By ``regular'' we mean that elements should be almost equilateral and equiangular. Regular elements have a positive effect on the stability of the system (see, for example, the numerical experiment in Section \ref{Sec:Regularity}). If an element is very irregular and very stretched in one direction, the approach of iterative subdivisions may fail. 
\subsubsection{Regular Quadrature}
If the element and the collocation node are sufficiently far apart from each other, \numcalc{} determines the order of the Gauss quadrature over the element by using an a-priori error estimator similar to \cite{LacWat76}. The number of quadrature nodes is determined by the smallest $m$ for which all of the 3 a-priori error estimates 
  \begin{align*}
  \varepsilon_G &:= 32\left(\frac{1}{2\tilde{r}}\right)^{2m+1},\\
  \varepsilon_H &:= 64(2m+1)\left(\frac{1}{2\tilde{r}}\right)^{2m+2},\\
  \varepsilon_E &:= 128(m+1)(2m+1)\left(\frac{1}{2\tilde{r}}\right)^{2m+3}
  \end{align*}
are smaller than the default tolerance of $10^{-3}$. If $m$ becomes larger than 6, $m$ is set to 6, which provides a compromise between accuracy and computational effort. Because of the splitting of elements into subelements done in the quasi-singular case, $m = 6$ is a feasible upper limit even in the regular quadrature.

For quadrilateral elements the quadrature nodes are given by the cross product of $m \times m$ Gauss nodes on the interval $[-1,1]$,~\cite[Table 25.4]{AbrSte64}. For triangular elements Gauss quadrature formulas for triangles based on $m$ are used ~\cite{Cowper73, Bathe02}.

\subsection{Linear System of Equations Solvers}
\numcalc{} offers to choose between a direct solver (either from LAPACK~\cite{Andersonetal99}, if \numcalc{} is compiled with LAPACK-support, or a slower direct LU-decomposition) or an iterative conjugate gradient squared (CGS) solver~\cite{Saad03}. The iterative solver is preferred because the number of elements $N$ is generally too big for a feasible computation with a direct solver. When using this solver, \numcalc{} offers to precondition the system matrix by applying either row-scaling or an incomplete LU-decomposition~\cite[Chapter 5]{Meister99}. However, the CGS solver may have robustness problems if the elements of the mesh are irregular. Such problems may result in convergence issues, even with preconditioning. When using the iterative solver:
\hint{Check if the iterative solver has converged.}
During the calculation \numcalc{} provides users with some information about the status of the iterative solver and the current residuum. Additionally, information about the iterative solver is stored in the output file  \texttt{NC.out} (see also \ref{App:Output}), e.g.,
\texttt{CGS solver: number of iterations = 25, relative error = 9.25255e-10.} If the relative residuum after a fixed number of iterations (in the default setting 250) is not below $10^{-9}$, a warning will be displayed and stored in \texttt{NC.out}. To help users of Mesh2HRTF, \texttt{NC.out} is automatically scanned for such warnings and users are explicitly warned.
\subsection{The Fast Multipole Method}\label{Sec:FMM_I}
One drawback of the BEM is that the system matrix is densely populated, which makes numerical computations with even moderately sized meshes prohibitively expensive in terms of memory consumption and computation time. This dense structure is caused by the Green's function $\frac{e^{\I k ||\bx -\by||}}{4\pi ||\bx - \by||}$ because the norm $||\bx - \by||$ introduces a nonlinear coupling between every element $\Gamma_j$  and every collocation node $\bx_i$.

The FMM is based on a man-in-the-middle principle where the Green's function $G(\bx,\by)$ is approximated by a product of three functions
\begin{equation}
G(\bx,\by) \approx G_1(\bx - \bz_1 ) G_2( \bz_1 - \bz_2) G_3( \bz_2 - \by).
\end{equation}
This splitting has the advantage that the integrals over the BEM elements $\Gamma_j$ can be calculated and used independently of the collocation nodes $\bx_i$. On the downside, this splitting is only numerically stable if $\bx$ and $\by$ are sufficiently apart from each other (for more details we refer to Sec.~\ref{Sec:Expansion}).

\begin{figure}[!h]
  \begin{center}
    \includegraphics[width=0.4\textwidth]{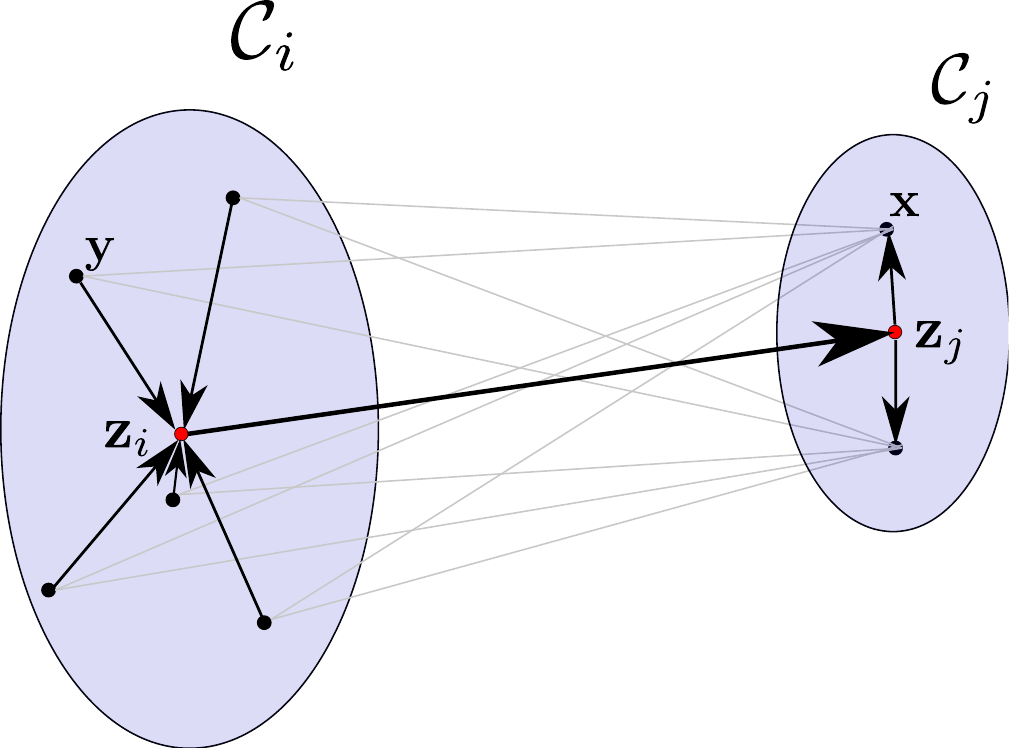}
  \end{center}
  \caption{Scheme of a FMM interaction between two clusters. The interactions (gray lines) between elements (black dots) located in two distant clusters (shaded areas) are reduced to fewer within and between cluster interactions  (black lines).}\label{Fig:FMMScheme}
\end{figure}

To derive this splitting, the elements of the mesh are grouped into different clusters $\Cl_i$. A cluster is a collection of elements that are within a certain distance of each other and for each cluster $\Cl_i$ its midpoint $\bz_i$ is defined as the average over the coordinates of all vertices of the elements in the cluster. For two clusters $\Cl_i$ and $\Cl_j$ with $\bz_i$ and $\bz_j$ that are sufficiently ``far'' apart, the element-to-element interactions of all the elements between both clusters are reduced to the local interactions between each elements of a cluster with its midpoint (functions $G_1$ and $G_3$), and the interactions between the two clusters (function $G_2$), respectively (cf. Fig.~\ref{Fig:FMMScheme}). Note that local-to-cluster and the cluster-to-cluster expansions can be calculated independently from each other.
\subsubsection{The expansion of the Green's function}\label{Sec:Expansion}
The expansion of the Green's function is based on its representation as a series~\cite{Rahola96}:
\begin{align}
  \nonumber G(\bx,\by) &= \frac{e^{\I k ||\bx - \by||}}{4\pi ||\bx - \by||} \\
  \nonumber &= \frac{\I k}{4\pi}\int\limits_\sphere e^{\I k (\bx - \bz_1 + \bz_2 - \by)\cdot \bs} \sum_{n = 0}^\infty \I^n (2n + 1) h^{(1)}_n(k||\bz_1 - \bz_2||) P_n \left( \frac{(\bz_1 - \bz_2)\cdot \bs}{||\bz_1 - \bz_2||}\right) \de\bs\\
  \label{Equ:FMMBase}&\approx \frac{\I k}{4\pi} \int\limits_\sphere  \Es(\bx,\bz_1,\bs) \D_L(\bz_1 - \bz_2,\bs) \T(\by,\bz_2,\bs) \de\bs,
\end{align}
where
\begin{align*}
  \Es(\bx,\bz_1,\bs) &= e^{\I k(\bx - \bz_1)\cdot \bs},\quad \T(\by,\bz_2,\bs) = e^{\I k (\bz_2 - \by)\cdot \bs},\\
  \D_L(\bz_1,\bz_2,\bs) &= \sum_{n = 0}^L \I^n (2n + 1) h_n^{(1)}(k||\bz_1 - \bz_2||) P_n \left( \frac{(\bz_1 - \bz_2)\cdot \bs}{||\bz_1 - \bz_2||}\right),
\end{align*}
$h^{(1)}_n$ denotes the spherical Hankel function of the first kind of order $n$, $P_n$ is given by the Legendre polynomial of order $n$, and $\sphere$ denotes the unit sphere in 3D. %and $\bs$ is a node on the unit-sphere $\sphere$.

The integral over $\sphere$ needs to be calculated numerically by discretizing the unit sphere $\sphere$   using $2L^2$ quadrature nodes $\bs_j, j = 1,\dots,2L^2$, where $L$ is the length of the multipole expansion in Eq.~(\ref{Equ:FMMBase}). The elevation angle $\theta \in [0,\pi]$ is discretized using $L$ Gaussian quadrature nodes. For the azimuth angle $\phi\in [0,2\pi]$ $2L$ equidistant quadrature nodes are used~\cite{Coifmanetal93}. 

It is not trivial to find an adequate truncation parameter $L$ for the sum in Eq.~(\ref{Equ:FMMBase}) and a lot of approaches have been proposed~\cite{Coifmanetal93,CecDar13}. On the one hand, $L$ needs be high enough to let the sum in Eq.~(\ref{Equ:FMMBase}) converge. On the other hand, a large $L$ implies Hankel functions with high order, which  become numerically unstable for small arguments. Thus, the FMM is not recommended if two clusters are close to each other.

In~\cite{Coifmanetal93}, a semi-heuristic formula for finding the optimal $L$ is given by
$$
L = k d_{max} + 5 \ln (k d_{max} + \pi)
$$
for single precision and 

$$
L = k d_{max} + 10 \ln (k d_{max} + \pi)
$$
for double precision, where $d_{\max}$ is the maximum inside-cluster distance $||(\bz_i - \by) + (\bx - \bz_j)||\le 2 r_{\max}$ over all clusters (see also Fig.~\ref{Fig:FMMScheme}) and $r_{\max}$ is the radius\footnote{The radius of a cluster is the biggest distance between the vertices in the cluster and the cluster midpoint.} of the largest cluster. As $L$ depends also on the wavenumber $k$, the order of the FMM is also frequency dependent. \numcalc{} uses a similar bound:
\begin{equation}\label{Equ:FMMTruncation}
L = \max(8,2 r_{\max} k + 1.8 \log_{10}(2 r_{\max} k +  \pi)) \approx \max(8,2 r_{\max} k + 0.7817 \ln(2 r_{\max} k +  \pi)). 
\end{equation} 
Our numerical experiments have shown that a lower bound for $L$ can enhance the stability and accuracy (especially at lower frequencies) and that a factor of 1.8 in Eq.~(\ref{Equ:FMMTruncation}) is a good compromise between stability, accuracy, and efficiency. Sec.~\ref{Sec:Benchmark} presents an example where the effect of this factor will be shown. 

In Eq.~(\ref{Equ:FMMBase}), the functions $\Es$ and $\T$ represent the local element-to-cluster expansion inside each cluster, whereas the function $\D_L$ represents the cluster-to-cluster interaction. One key element in $\D_L$ are the spherical Hankel functions $h^{(1)}_n$ of order $n = 0,\dots, L$, that become singular at 0, $\lim_{x\rightarrow 0} h_n^{(1)}(x) = O\left(\frac{1}{x^{n+1}}\right)$~\cite[Eq. 10.1.5]{AbrSte64}, see also Sec.~\ref{Sec:FMM-Matrices}.  The FMM expansion is therefor only stable for relatively large arguments $k||\bz_1 - \bz_2||$, which implies that the FMM can only be applied for cluster-pairs that are sufficiently apart from each other. Following~\cite{CecDar13}, \numcalc{} analyzes the relation between the distance $||\bz_i - \bz_j||$ of the cluster midpoints to their radii $r_i$ and $r_j$, which are given by the maximum distance between the cluster midpoint and the vertices of the cluster elements. If
\begin{equation}\label{Equ:Farfieldcond}
||\bz_i - \bz_j|| > \frac{2}{\sqrt{3}} (r_i + r_j),
\end{equation}
the two clusters are defined to be in each others farfield and the FMM expansion is applied. 

If two clusters do not fulfil Eq.~(\ref{Equ:Farfieldcond}), the clusters and their elements are defined to be in each others nearfield. The interaction between such cluster pairs will be calculated using the conventional BEM approach with the non-separable $G(\bx,\by)$ leading to the sparse nearfield matrix $\bN$. If two clusters are found to be in each others farfield, the Green's function $G(\bx,\by)$ will be approximated by Eq.~(\ref{Equ:FMMBase}).

Note that the multipole expansions for the derivatives of $G$ are easy because $\bx$ and $\by$ only occur as arguments of  exponential functions, for which the calculation of the derivative is trivial.
\subsubsection{Cluster generation}\label{Sec:SLFMM}
To actually generate the clusters in \numcalc{}, a bounding box is put around the mesh. This box is then subdivided into (approximately) equally sized sub-boxes for which the average edge length can be provided by the user. Because the sub-boxes have the same edge length in all three dimensions, their number may vary along the different coordinate axes.

If the edge length is not specified, the default edge length will be $|e_\text{b}| = \left(\sqrt{N}A_0\right)^{1/2}$, where $N$ is the number of elements and $A_0$ is the average area of all elements. The number of subdivisions is estimated by rounding the quotient of the edge length of the original box and the target edge length. The actual edge lengths of the sub-boxes can slightly deviate from the target length $|e_\text{b}|$. If the midpoint of an element is within a sub-box, the element is assigned to the box, and all elements inside the $i$-th sub-box form the cluster element $\mathcal{C}_i$. The choice of initial edge length can be motivated by the following idea: First, it is assumed that all elements have roughly the same size. Secondly, the assumption of having about $n_0 = \sqrt{N}$ clusters with about $N_0 = \sqrt{N}$ elements in each cluster provides a nice balance between the number of the clusters and the number of elements in each cluster. If we assume that the surface of the scatterer is locally smooth, the elements in each cluster cover approximately a surface area of $A_\Cl = N_0 A_0$. As each cluster is fully contained in one bounding box, a good estimate for the edge length of such a box is $\sqrt{A_\Cl} = \sqrt{N_0 A_0} = \sqrt[4]{N} \sqrt{A_0}$. Note that the assumption of having elements of roughly the same size has an impact on the efficiency of the method, but also influences the stability of the FMM expansion. If the size of the elements inside a cluster varies too much, the local element-to-cluster expansions may run into numerical troubles. This means that
\hint{The elements should have approximately the same size, at least locally.}
Sub-boxes that contain no elements are discarded. Clusters with a small radius (in comparison with the average cluster radius) are merged with the nearest cluster with large size.  
\begin{figure}[!h]
  \begin{center}
    \includegraphics[width=0.3\textwidth]{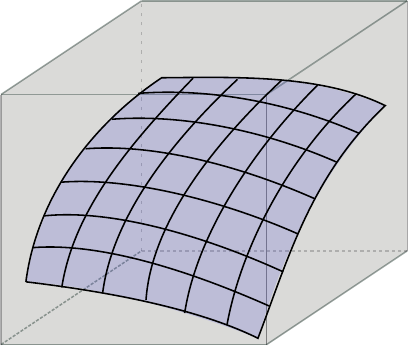}
    \caption{Example of a cluster (blue) and the bounding (sub)box (grey) around it.}\label{Fig:ClustervsBox}
  \end{center}
\end{figure}

Note that there is a logical difference between the sub-box and the cluster. A cluster is just a collection of elements, which form a 2D manifold, whereas the bounding box is a full 3D object (see also Fig.~\ref{Fig:ClustervsBox}). %, which is just a collection of the elements.
As a consequence, the radius of the cluster is usually \emph{not} the same as half the edge length of the sub-box. Also, the cluster midpoint is usually different from the midpoint of the sub-box. 
\subsubsection{FMM System of Equations}
Combining the FMM expansion Eq.~(\ref{Equ:FMMBase}) with the discretized BIE  Eq.~(\ref{Equ:SystemI}), the system matrix for the FMM has the form 
\begin{equation}\label{Equ:SystemSLFMM}
\bN + {\bf S D} {\bf T},
\end{equation}
in which $\bN$ contains the contributions for cluster-pairs in the near-field, where the FMM cannot be applied, and ${\bf S, D}$ and ${\bf T}$ contain the contributions of far-field cluster pairs, where the FMM expansion will be applied. Because of the clustering, the matrices ${\bf S, D}$ and ${\bf T}$ have a block-structure, see Fig.~\ref{Fig:Blockstructure} for a schematic  representation of block structure of ${\bf SDT}$.

To better illustrate the sparse structure of these matrices, we look at the entries of the different matrices with respect to the integral parts involving the Green's function $G(\bx,\by$), the matrix parts for the integrals with respect to $H, H',$ and $E$ can be derived in a similar way.

The entries of the near-field matrix $\bN$ are calculated for each collocation node $\bx_i$ using the quadrature techniques for standard BEM (see Sec.~\ref{Sec:Quadrature}), resulting in
  $$
  \bN_{ij} = \left\{
    \begin{array}{cc}
      \int\limits_{\Gamma_j} G(\bx_i,\by)\de\by & \Gamma_j \in \N_i,\\
      0 & \text{otherwise},
    \end{array}
  \right. 
  $$
where $\Gamma_j$ is in the nearfield $\N_i$ of the cluster containing the collocation node, and $i,j = 1,\dots,N$. 

The matrix ${\bf D}$ represents the cluster-to-cluster interaction between Clusters $\Cl_m$ and $\Cl_n$. It contains $2L^2$ blocks of size $N_C \times N_C$ defined by
  $$
  (\bD_\nu)_{mn} := \left\{
    \begin{array}{cl}
      \D_L(\bz_m,\bz_n,\bs_\nu) &\text{ if } (C_m,\Cl_n) \text{ is a far-field cluster pair},\\
      0 &\text{ otherwise},
    \end{array}
    \right. 
  $$
where $N_C$ is the number of clusters, $m,n = 1,\dots N_C$, and $\nu = 1,\dots,2L^2$, with $2L^2$ being the number of quadrature nodes on the unit sphere. The total size of ${\bf D}$ is $2 N_C L^2 \times 2N_CL^2$. 

The matrix $\bT$ describes the local element-to-cluster expansion from each element of a cluster to the cluster midpoint and is represented as a ${2L^2}\times 1$ dimensional block matrix
  $$
  \bT = (\bT_1,\cdots,\bT_{2L^2})^\top,
  $$
  where each single block has size $N_C \times N$ and is defined by
  $$
  (\bT_\nu)_{mj} = \left\{
    \begin{array}{cc}
      \int\limits_{\Gamma_j}e^{\I k(\bz_m - \by) \cdot \bs_\nu}\de\by & \Gamma_j\in \Cl_m, \\
      0 & \text{otherwise}.
    \end{array}
    \right. 
  $$
The local distribution matrix $\bS$ of size $1\times 2L^2$ is defined by 
  $$
  \bS = (\bS_1,\dots,\bS_{2L^2}),
  $$
where each $N\times N_C$ block has the form
  $$
  (\bS_\nu)_{in} = \left\{
    \begin{array}{cc}
      w_\nu e^{\I k(\bx_i - \bz_n)\cdot \bs_\nu} & \bx_i \in \Cl_n,\\
      0 & \text{otherwise},
    \end{array}
  \right.
  $$
where $w_\nu$ is the $\nu$-th weight used by the quadrature method for the integral over the unit sphere.
\begin{figure}
  $$
  \left(
\begin{array}{ccc}
  \fbox{\parbox[][100pt][c]{20pt}{%
  \phantom{ai}$\bS_1$
  }
  } &
  \cdots &
  \fbox{\parbox[][100pt][c]{20pt}{%
  \phantom{i}$\bS_{2L^2}$
  }
  }
\end{array}
\right)
\left(
  \begin{array}{ccc}
    \fbox{\parbox[][20pt][c]{20pt}{%
    \phantom{a}$\bD_1$
    }}
    & 0 & 0\\
    0 & \ddots & 0 \\
    0 & 0 &
            \fbox{\parbox[][20pt][c]{20pt}{%
            $\bD_{2L^2}$
            }}
  \end{array}
\right)
\left(
  \begin{array}{c}
    \fbox{\parbox[][20pt][c]{100pt}{%
    \centering
    $\bT_1$
    }
    }
    \\[13pt]
    \vdots\\[13pt]
    \fbox{\parbox[][20pt][c]{100pt}{%
    \centering
    $\bT_{2L^2}$
    }
    }
  \end{array}
\right)
$$
\caption{Block structure of the FMM matrices in the SLFMM. Each matrix block $\bS_i$ has size $N\times N_C$, where $N$ is the number of elements and $N_C$ is the number of clusters. Each block $\bD_i$ has size $N_C\times N_C$ and each block $\bT_i$ has size $N_C \times N$.}\label{Fig:Blockstructure}
\end{figure}
\subsubsection{Multi level FMM (MLFMM)}
For higher efficiency, the FMM can be used in a multilevel form, in which a cluster tree is generated. The mesh is not only subdivided once in the first level (root level), but clusters are iteratively subdivided in the next levels until the clusters at the finest level (leaf level) have about 25 elements. As the FMM is based on multiple levels, we speak of the multilevel fast multipole method (MLFMM) compared to the single level fast multipole method (SLFMM) described in Sec.~\ref{Sec:SLFMM}.

In contrast to other implementations of the MLFMM where the cluster on the root level contains all elements of the mesh, the root of the cluster tree in \numcalc{} is based on a clustering similar to the one used for the SLFMM (see Sec.~\ref{Sec:SLFMM}). Users can either define an initial bounding box edge length $|e_{\text{r}}|$ or use the default target edge length chosen such that each cluster contains approximately $N_0 = 0.9 \sqrt{N}$ elements. As in the SLFMM, it is assumed that the elements have similar sizes with an area of $A_0$ and an average edge length of $|e_0|$.
\begin{figure}[!h]
  \begin{center}
    \includegraphics[width=0.25\textwidth]{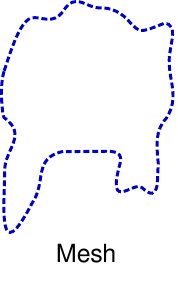}
    \hspace{10pt}
    \includegraphics[width=0.25\textwidth]{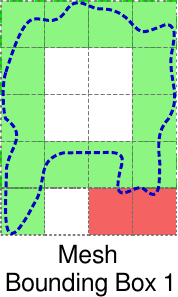}
    \hspace{10pt}
    \includegraphics[width=0.25\textwidth]{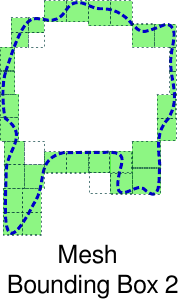}
    \caption{Example of the MLFMM clustering: Mesh (dark/blue line) and the bounding boxes for the first two levels. At each level, empty (white) boxes are discarded, small (dark grey/red) boxes are merged into the neighbouring large (light grey/green) boxes, which are used for the clustering. }\label{Fig:Clustering}
  \end{center}
\end{figure}

The clusters at higher levels are constructed in two steps: First, new bounding boxes are created for each cluster based on the coordinates of element vertices in the respective cluster. Then, sub-boxes with target edge length $\frac{|e_{\text{r}}|}{2^\ell}$ are created similar to the root level, where $\ell$ is the number of the level. The elements within these sub-boxes then form the clusters in the new level of the MLFMM. The clusters contained in the sub-boxes in the new level are called children, the original cluster element is called parent cluster. Fig.~\ref{Fig:Clustering} shows an example of creating the bounding boxes and clusters for the root and next level.

Currently, \numcalc{} does \emph{not} build the cluster tree adaptively, but the number of levels
$$
\ell_{\max} = \max\left(1, \text{round}\left[\log_2\left(\frac{|e_{\text{r}}|}{|e_{\text{l}}|}\right)\right]\right) = \max\left(\text{round}\left[\frac12 \log_2\left(\frac{\sqrt{N}}{22.5}\right)\right]\right)
$$ 
is determined beforehand, where $|e_{\text{r}}|$ and $|e_{\text{l}}|$ are the target bounding box lengths on the root and leaf level, respectively. Because in the default setting each cluster on the root level contains $N_0 = 0.9 \sqrt{N}$ elements, the edge length of each bounding box can be estimated by $|e_{\text{r}}| = \sqrt{\frac{N}{N_0}} \sqrt{A}$. On the leaf level each bounding box contains about 25 elements and its edge length is roughly $|e_{\text{l}}| = 5 \sqrt{A}$. In that setting, $\ell_{\max}$ defines how often $|e_{\text{r}}|$ needs to be split to obtain $|e_{\text{l}}|$.  \numcalc{} produces cluster trees which are well-balanced and where all leaf clusters are on the same level. Note that the bounding boxes are constructed to have about the same edge length in all directions. If the mesh has large aspect ratios, the number of clusters in the respective dimension will differ.

In practice, elements are rarely distributed regularly inside the bounding boxes. Thus, our clustering approach may lead to leaf clusters with very few elements on the one hand and clusters with more then 25 elements on the other hand. Our numerical experiments have shown that having already about $\sqrt{N}$ clusters in the root level improves the stability and efficiency of the FMM, because, for example, the expansion length used by the FMM  (see Sec.~\ref{Sec:Expansion}) is kept relatively low.

In order to solve the MLFMM case, nearfield and farfield cluster-pairs are defined on each level using the same rule as for the SLFMM. However, in contrast to the SLFMM, the FMM is \emph{not} applied for all far-field cluster pairs on each level, but only to a specific sub-set, the so called interaction list. For each cluster $\Cl_j$ the interaction list $I(\Cl_j)$ is defined as the elements in the child clusters of the near field clusters of the parent of $\Cl_j$.

\begin{figure}[!h]
  \begin{center}
    \includegraphics[width=0.8\textwidth]{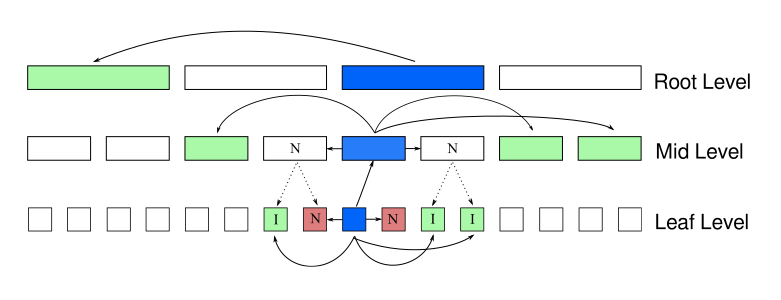}
  \end{center}
  \caption{Scheme of the MLFMM. For a single cluster on the highest level $\ell_{\max}$ (= leaf level, shown in the bottom row by the dark blue box), the red boxes denote its nearfield clusters and the light green boxes denote the clusters stored in the interaction list, for which the FMM is applied. The interaction with the other clusters within this level (denoted by the white boxes) are considered in the parent level $\ell_{\max}\!-\!1 $. This procedure is repeated for all levels of the MLFMM. Arrows depict interactions between clusters and messages passing to the parent, dotted arrows depict a parent-child relation.}\label{Fig:MLFMM}
\end{figure}

Calculations for the MLFMM start at the leaf level $\ell_{\max}$. For each  cluster $\Cl_i^{\ell_{\max}}$ on this level the near field matrix ${\bf N}$ is calculated for all near field cluster pairs as in the SLFMM case (for example, the blue box and the two red boxes in Fig.~\ref{Fig:MLFMM}). The fast multipole expansion, however, is only applied between each cluster $\Cl^{\ell_{\max}}_i$ and the clusters in the respective interaction list $I(\Cl_i^{\ell_{\max}})$ (green boxes in Fig.~\ref{Fig:MLFMM} for the selected blue cluster box). All other cluster-pairings (white boxes) on this level are neglected.

For those pairs the local element-to-cluster expansions at level $\ell_{\max}$ are transformed into local element-to-cluster expansions with respect to the parent cluster in level $\ell_{\max}\!-\!1$ (upward pass, blue box to parent blue box in Fig.~\ref{Fig:MLFMM}). This transformation is necessary because the parent cluster has a different cluster midpoint and maybe a different multipole expansion length as the child. At this level no near-field components need to be calculated, because these have already been ``covered'' by the calculations at the level of the children. For each pairing of cluster $\Cl^{\ell_{\max}-1}_i$ and clusters in its interaction list, the FMM is performed. For all other cluster pairs the local expansions are again passed on to the parent, and the procedure is repeated at the level above.

In summary, at each level $\ell$ the local element-to-cluster matrices $\bT(\ell)$, the local cluster-to-element matrices $\bS(\ell$), and the cluster-to-cluster interaction matrices $\bD(\ell)$ will be calculated, and the final system for the MLFMM is given by
$$
N(\ell_{\max}) + \sum_{\ell = 1}^{\ell_{\max}} \bS(\ell) \bD(\ell) \bT(\ell).
$$

\begin{figure}[!h]
  \begin{center}
    \includegraphics[width=0.5\textwidth]{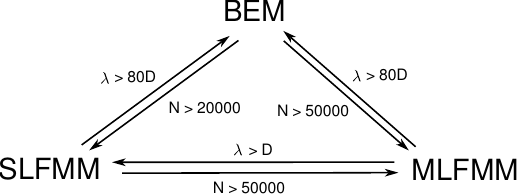}
  \end{center}
  \caption{Scheme of the automatic changes in the FMM-versions. Here $\lambda$ is the wavelength and $D$ is the diameter of the object, i.e., the biggest distance between two points of the mesh.}\label{Fig:FMMChange}
\end{figure}

While the MLFMM is recommended for most of the problems, \numcalc{} implements a few ``safeguards'' focused on an efficient calculation and stable results (see also Fig.~\ref{Fig:FMMChange}):
\begin{itemize}[leftmargin=*]
\item If the user selects the conventional BEM (i.e., no FMM) but the number of elements $N$ is larger than $20000$, \numcalc{} will switch to the SLFMM.
\item If the user selects the conventional BEM or the SLFMM and the number of elements $N$ is larger than $50000$, \numcalc{} will automatically switch to the MLFMM.
\item If the user selects the FMM, but the wavelength $\lambda$ is $80$ times larger (or more) than the diameter $D$ of the scatterer (the maximum distance $r_{\max}$ between two nodes of the mesh), \numcalc{} will use the conventional BEM (i.e., \emph{without} the FMM). This may especially the case at low frequencies and prevents the too small truncation of the multipole method (see Section~\ref{Sec:Expansion}).
\item If the user selects the MLFMM, but $\lambda$ is larger than $r_{\max}$, \numcalc{} automatically switches back to the SLFMM. 
\end{itemize}
The current version of \numcalc{} triggers no explicit warning with respect to that switch, however, in the output file \texttt{NC.out} information about the method used for each frequency step is provided. 

One specific detail in \numcalc{} is the fact that \numcalc{} determines the matrices $\bT(\ell), \bD(\ell)$ and $\bS(\ell)$  beforehand for \emph{all} levels and explicitly keeps those sparse matrices in memory.
This is in contrast to many other implementations of the FMM that only calculate the expansions on the leaf level and then pass the information between levels by interpolation and filtering algorithms (cf.~ \cite{AmiPro03}).
The approach in \numcalc{} has the advantage of a faster calculation because no interpolation and filtering routines are necessary. 

Nevertheless, one drawback of storing $\bT(\ell)$ and $\bS(\ell)$ on each level $\ell \le \ell_{\max}$ is that these matrices may become large close to the root level, which can cause problems with memory consumption. The radii of the clusters are larger at smaller levels (close to the root of the cluster tree). %grows for smaller levels $\ell$,
This means that the truncation parameter $L$  in the multipole expansion (see Eq.~(\ref{Equ:FMMTruncation})) and the number of necessary quadrature nodes for the integral over the unit sphere increase when getting closer to the root level. Thus, the number of blocks in $\bT(\ell)$ and $\bS(\ell)$ also increases towards the root, which means higher memory consumption. However, with the default clustering of having $\sqrt{N}$ clusters in the root level, large clusters can be avoided  and the memory consumption can be kept at a feasible level. Only for very large meshes with more then 100.000 to 150.000 elements the additional memory requirement can become quite large (see, e.g., the problem presented in Sec.~\ref{Sec:MLFMM4}).  
\section{Critical components of \numcalc{} and Their Effect}\label{Sec:Critical}
While users may not be necessarily interested in the details covered in Sections~\ref{Sec:BriefBEM} and \ref{Sec:Implementation}, they should be aware of some of the consequences of the assumptions used in these sections and the key components of the BEM and the FMM.

\subsection{Singularity of the Green's function}\label{Sec:GreenSgl}
The Green's function $G(\bx,\by)$ and its derivatives $H(\bx,\by), H'(\bx,\by),$ and $E(\bx,\by)$ are essential parts of the BEM. These functions become singular whenever $\bx = \by$, i.e., when we integrate over the element $\Gamma_i$ that contains the collocation point $\bx_i$. There are algorithms that can deal with these type of singular integrals (cf.~\cite{Duffy82, Krishnasamyetal90}), but the singularities also cause \emph{numerical} problems for elements in the neighbourhood of $\Gamma_i$. In \numcalc{}, these almost-singular integrals are calculated by subdividing the element several times (see also Section \ref{Sec:Quadrature}), but these subdivisions slow down calculations. The (almost) singularity especially causes problems if
\begin{itemize}%[leftmargin=*]
\item an external sound source is close to the surface of the object,
\item there are very thin structures in the geometry, thus, if front and backside of the object are close to each other,
\item there are overlapping or twisted elements.
\end{itemize}
If \numcalc{} terminates with an error code message ``\texttt{number of subels.\ which are subdivided in a loop must <= 15}'', the reason is most likely the presence of irregular or overlapping elements causing more then 15 subdivisions of an element. \numcalc{} also displays the element index for which the error occurred. There is a high chance that there is some problem with the mesh or an external sound source close to this element. So:
\hint{If you get a subdivision loop error, check your mesh in the vicinity of the element for overlaps and irregularities. If there is a external sound source close to that element, move the source further away from the mesh.}
\subsection{FMM matrices}\label{Sec:FMM-Matrices}
The expansion Eq.~(\ref{Equ:FMMBase}) of the Green's function has several consequences especially in connection with finding the right truncation length, see Section~\ref{Sec:Expansion}.
For small arguments, the spherical Hankel function $h^{(1)}_n(x) = j_n(x) + \I y_n(x)$ of the first kind and of order $n$ behaves like~\cite[Eqs.~(10.1.4),(10.1.5)]{AbrSte64}
$$
\lim_{x\rightarrow 0} j_n(x) = \frac{1}{1\cdot3\cdot5\cdots (2n+1)}{x^n},\quad
\lim_{x\rightarrow 0} y_n(x) = -\frac{1\cdot3\cdot5\cdots(2n-1)}{x^{n+1}}.
$$
Numerical problems may arise for large expansion lengths and small arguments because the spherical Hankel function goes to infinity with $O(\frac{1}{x^{n+1}})$. This is a problem when $k||z_1 - z_2|| \ll 1$, which can happen at low frequencies and clusters being close to each other.

The spherical Hankel functions are implemented with a simple recursive procedure~\cite[Section 3.2.1]{Giebermann97} with the aim to find a balance between accuracy and computation time. For higher orders the calculations become more involved and, as the calculation is based on a recursive formula, errors from lower orders propagate and add up at higher orders. Hence, the length of the FMM expansion has a crucial effect on the results of the FMM:
\hint{If you need an expansion order $L$ greater 30, check the convergence of the iterative solver and your results for plausibility.}
Like the number of iterations the expansion length $L$ is given in \texttt{NC.out} for each level of the FMM.

The truncation parameter $L$ also directly influences the number of quadrature nodes on the sphere (see Sec.~\ref{Sec:Expansion}) and the memory required to store the local-to-cluster expansion matrices $\bT$ and $\bS$. Because $L$ depends on the cluster radius (see  Eq.~\ref{Equ:FMMTruncation}) one remedy for problems caused by a large $L$ is to adjust the length of the initial bounding box such that the radii of the clusters become smaller. 

\section{Benchmarks}\label{Sec:Benchmark}
We present results for three benchmark problems in acoustics to illustrate the limitations of \numcalc{} in terms of accuracy and computer resources. First, we analyze the scattering of a plane wave on a sound-hard sphere. Second, we analyze the calculation of the sound pressure on a human head caused by point-sources placed around the head. Finally, we analyze the calculation of the sound field inside a duct. 

\subsection{Sound-Hard Sphere}
We consider the acoustic scattering of a plane wave on a sound-hard sphere with radius $r = 1$\,m at different frequencies $f = 200, 800,$ and $1000$\,Hz. The speed of sound is assumed to be $c = 340$\,m/s,  the wavenumber $k \approx 3.696, 14.784, 18.48$\,rad/m, and the wavelength $\lambda = 1.7, 0.425$, and $0.34$\,m. For this problem an analytic solution based on spherical harmonics is known~\cite[Eq.~(6.186)]{Williams99}. To investigate the effect of various discretization methods, the meshes of the sphere were created using two methods. The first type of meshes (cube-based meshes) was created by discretizing the unit cube with triangles and by projecting the discretized cube onto the unit sphere. The second type of mesh (icosahedron-based meshes) was the often-used construction using an (subdivided) icosahedron that is projected onto the sphere~\cite{Ward15}. Both meshes were created using Matlab/Octave scripts. 

\begin{figure}[!h]
  \begin{center}
    \begin{tabular}{rcrc}
      \raisebox{0.35\textwidth}{a) \hspace{-30pt}} &
      \includegraphics[width=0.48\textwidth]{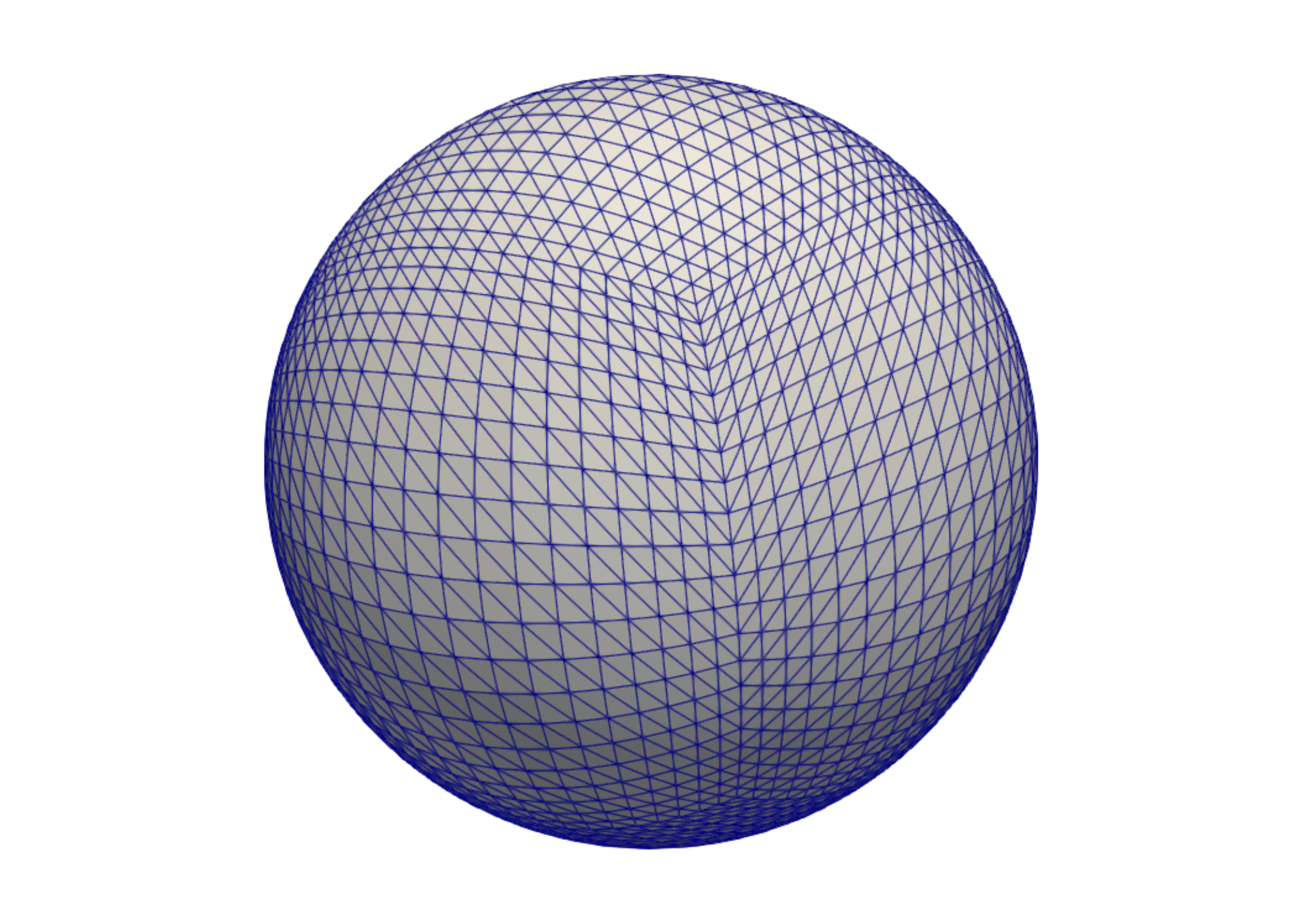} &
      \raisebox{0.35\textwidth}{b) \hspace{-30pt}} &
      \includegraphics[width=0.48\textwidth]{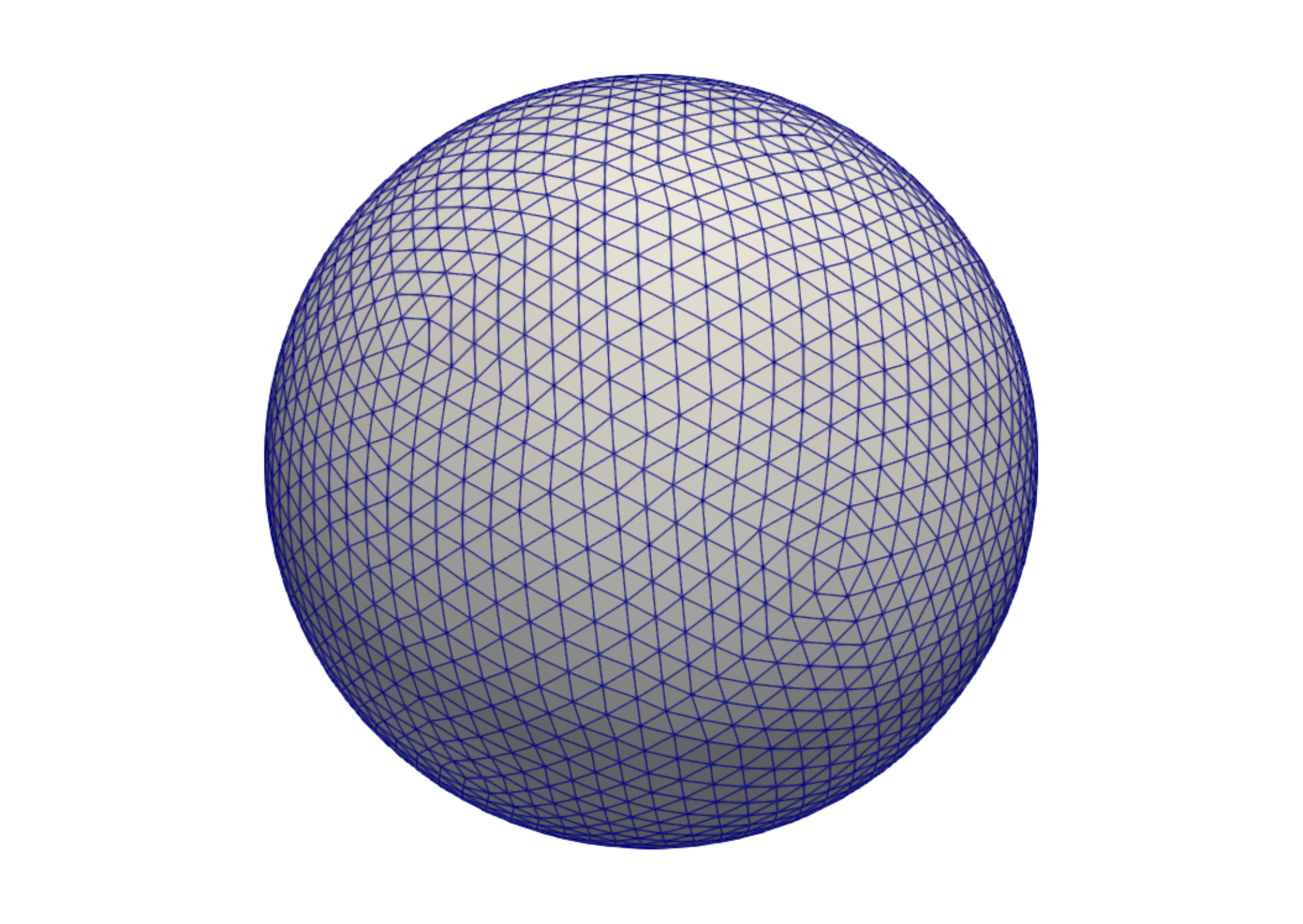}
    \end{tabular}                                                    
    \caption{Discretizations of the unit sphere used in our benchmark. On the left side, the mesh was generated by discretizing a cube and projecting it onto a sphere. Note the smaller triangles close to the projected edges of the cube. On the right, the mesh is based on a subdivided and projected icosahedron.}\label{Fig:Spheremesh}
  \end{center}
\end{figure}

In general, the icosahedron-based meshes had almost equilateral triangular elements of about the same size, but the number of elements were restricted to $N_{\text{ico}} = 20\cdot4^n, n\in \mathbb{N}_0$. The number of elements for the cube-based meshes were only restricted to $N_{\text{cube}} = 6\cdot (n+1)^2, n\in \mathbb{N}_0$, but the triangles were not as regular as in the icosahedron-based meshes and the size of the triangles differed over the sphere, especially near the projected edges of the cube. An example for both meshes is depicted in Fig.~\ref{Fig:Spheremesh}, where the sphere on the left side is based on a projected cube and has 5292 elements. The mesh on the right is based on a icosahedron and has 5120 elements. Both meshes fulfill the 6-to-8 elements per wavelength rule  for frequencies smaller than roughly 950\,Hz. 

\numcalc{} calculates the acoustic field on the surface on the \emph{collocation} nodes, which, in this example, do not lie on the surface of the unit sphere. The errors depicted below are thus a combination of geometrical and numerical error. However, numerical experiments with a modified radius showed that the error caused by the geometry is negligible compared to the numerical error. 

\begin{figure}[!h]
  \begin{center}
  \includegraphics[width=0.48\textwidth]{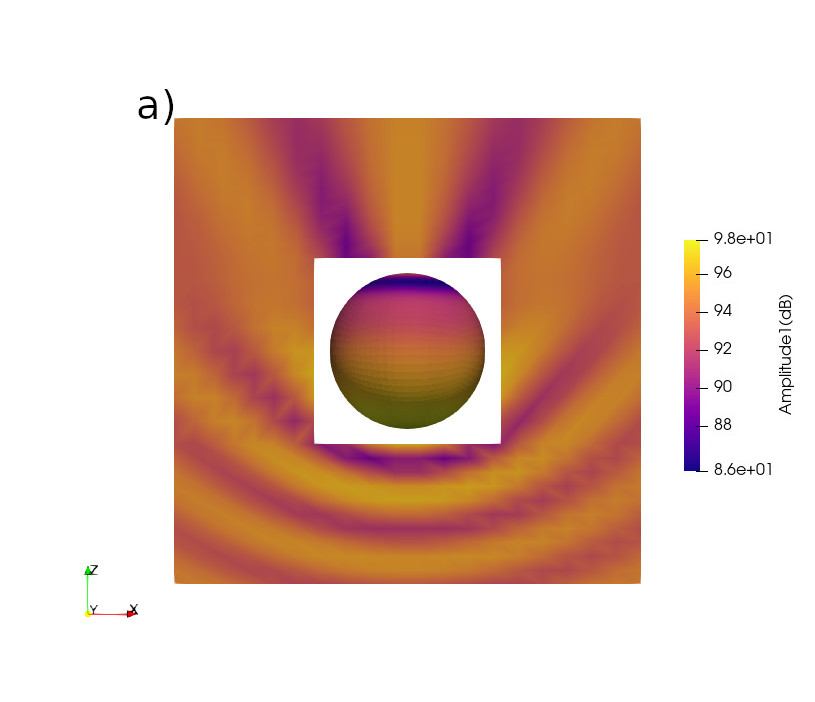}
  \includegraphics[width=0.48\textwidth]{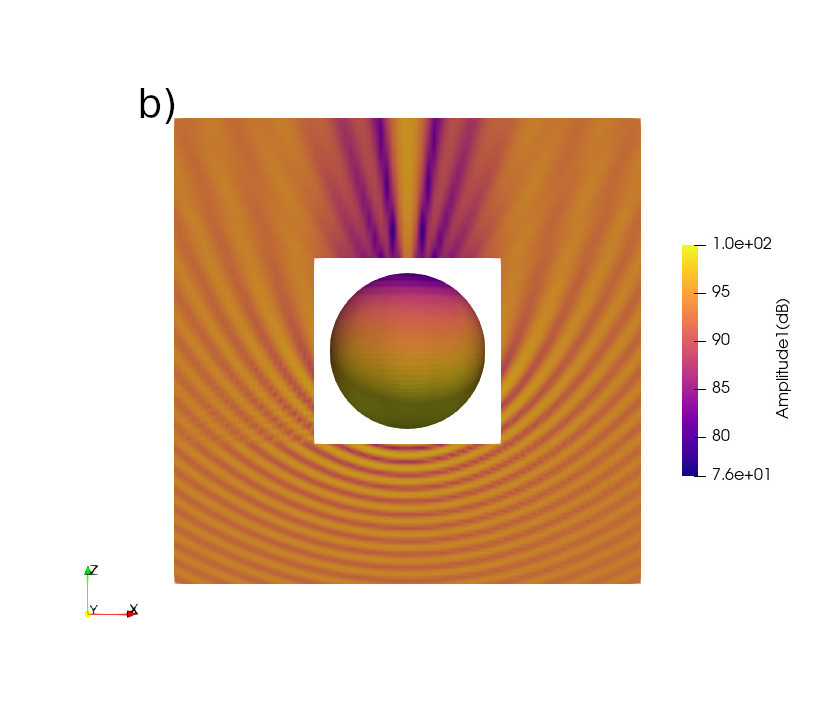}
  \caption{Example of the acoustic field calculated for 200\,Hz (left panel) and 1000\,Hz (right panel) caused by a plane wave propagating in the positive $z$-direction.}\label{Fig:BenchmarkField}
  \end{center}
\end{figure}

To illustrate the combination of incoming and scattered  acoustic field and to investigate the errors in the acoustic field outside the sphere, 24874 evaluation points have been placed on a plane around the sphere. Fig.~\ref{Fig:BenchmarkField} shows the logarithmic sound pressure level (SPL) in dB as an example of such acoustic field calculated for two frequencies. These evaluation points were also used to compute the errors between the BEM solution and the analytic solution outside the sphere. The MLFMM with default settings was used to produce the results for these figures.
\subsubsection{Effect of element regularity}\label{Sec:Regularity}
\begin{figure}[!h]
  \begin{center}
    \includegraphics[width=0.75\textwidth]{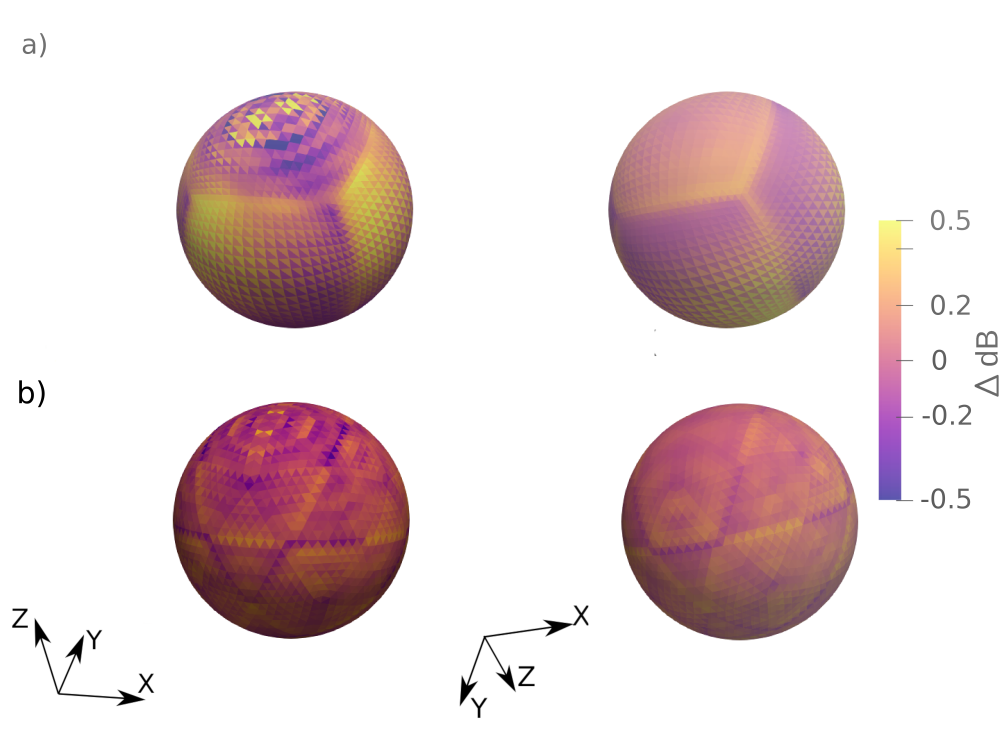}
    \caption{Sound pressure level differences (in dB) between the BEM calculation and the analytical solution at $f = 800$\,Hz. a) Cube-based mesh, b) Icosahedron-based mesh. As external source a plane wave in the negative $z$ direction was used.}\label{Fig:ErrorSphere}
  \end{center}
\end{figure}

Fig.~\ref{Fig:ErrorSphere} shows the difference between the numerical and the analytical solutions to illustrate the effect of element regularity. The color represents the SPL-difference between the analytic and the numerical solution at the midpoint of each element on the sphere. The effect is shown for the frequency of $f = 800$\,Hz and the two different types of discretizations. Both meshes are fine enough to fulfil the 6-to-8 elements per wavelength rule.

Despite having a similar number of elements, the different regularity in the two meshes has a clear effect on the errors. For the very regular icosahedron-based triangularization, the errors between the calculated and the analytic SPLs were between $\pm 0.34$\,dB. For the cube-based triangularization, the errors were larger, i.e., between $-0.7$ and $0.9$\,dB, being considerably larger at the triangles around the projected cube edges. For the cube-based mesh, we observed a clearer distinction in the errors between the ``sunny'' and ``shadow'' side, i.e., the side oriented towards the plane wave and the side oriented away from the plane wave. For the icosahedron-based mesh the error distribution was rather smooth across the two sides. 
\subsubsection{Error as function of the number of elements $N$}\label{Sec:Errororder}
In this section, we analyzed how the number of elements affect the relative error
$$
\varepsilon_r(\bx) := \frac{|p(\bx) - p_{0}(\bx)|}{|p_{0}(\bx)|}
$$
between the calculated acoustic field $p$ and the analytical solution $p_{0}$ given in \cite[Eq.~(6.186)]{Williams99}. We considered  the error on the surface as well as outside the sphere as a function of the number of elements $N$. Further, we used only cube-based meshes because these meshes offer more flexibility in terms of possible numbers of elements $N$.

\begin{figure}[!h]
  \begin{center}
  \begin{tabular}{rcrc}
      \raisebox{0.35\textwidth}{a)\hspace{-20pt}} &
      \includegraphics[width=0.46\textwidth]{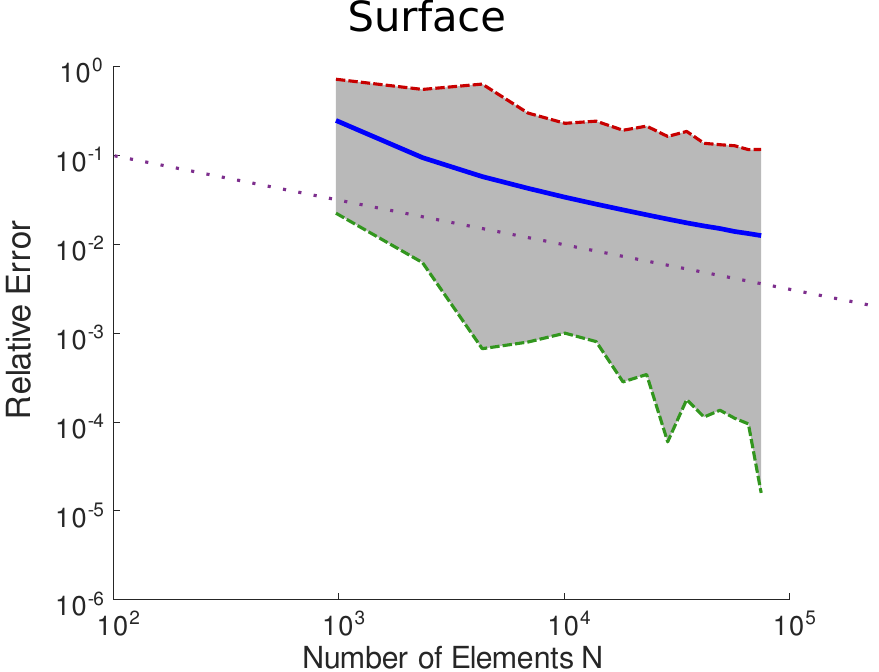} &
      \raisebox{0.35\textwidth}{b)\hspace{-20pt}} &
      \includegraphics[width=0.46\textwidth]{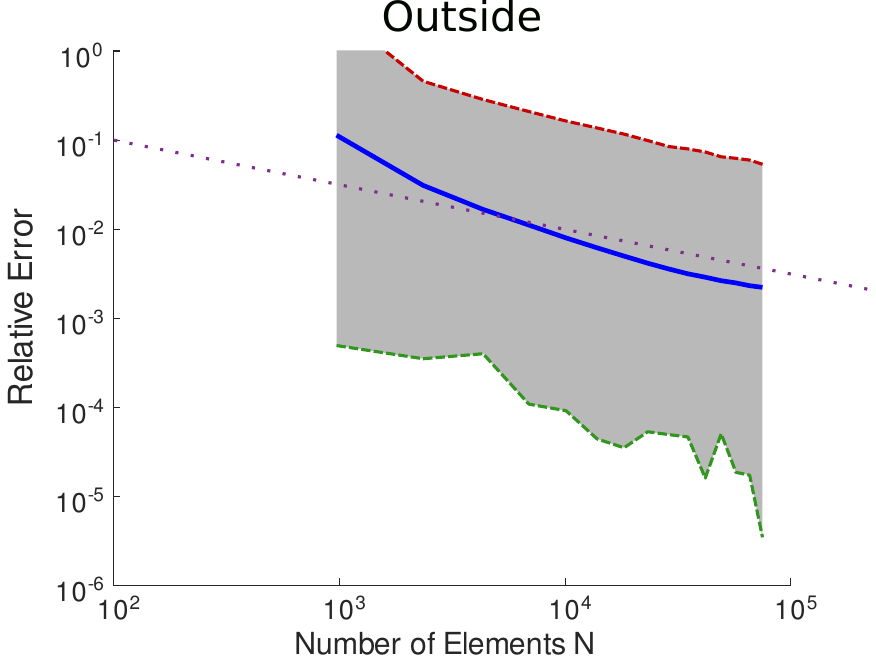}
  \end{tabular}
  \end{center}                                                       
  \caption{Relative errors as a function of the number of elements $N$ for the sphere benchmark at $f = 1000$\,Hz. Left side: Error on the surface of the sphere. Right side: Error around the sphere. The continuous (blue) line depicts the mean error over all nodes, the gray area depicts the area between the maximum and the minimal  error (dashed lines). The dotted line depicts the function $\frac{1}{\sqrt{N}}$.}\label{Fig:Error1000}
\end{figure}

Figs.~\ref{Fig:Error200} and~\ref{Fig:Error1000} show the maximum (red dashed line), mean (blue continuous line), and minimum (green dashed line) of the relative errors calculated at 200 and 1000\,Hz, respectively. It is known that the numerical error of the collocation method with constant elements is linearly proportional to the average edge length of the mesh. To illustrate this theoretical behaviour, the dotted line represents the functions $\frac{1}{\sqrt{N}}$, which is roughly proportional to the average edge length of the elements. The errors are in principle in the range of $O(\frac{1}{\sqrt{N}})$, however, they increased with the number of elements in the lower frequency example (Fig.~\ref{Fig:Error200}). For \numcalc{}, it can be shown that this effect is caused by the following specific details of the implementation:
\begin{figure}[!h]
  \begin{center}
    \begin{tabular}{rcrc}
      \raisebox{0.35\textwidth}{a) \hspace{-20pt}} &
      \includegraphics[width=0.46\textwidth]{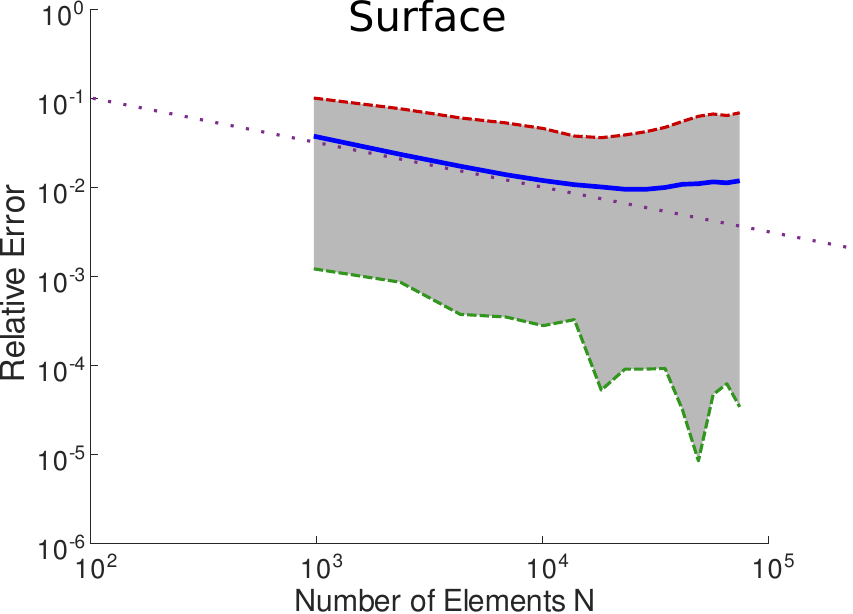} &
      \raisebox{0.35\textwidth}{b) \hspace{-18pt}} &
       \includegraphics[width=0.46\textwidth]{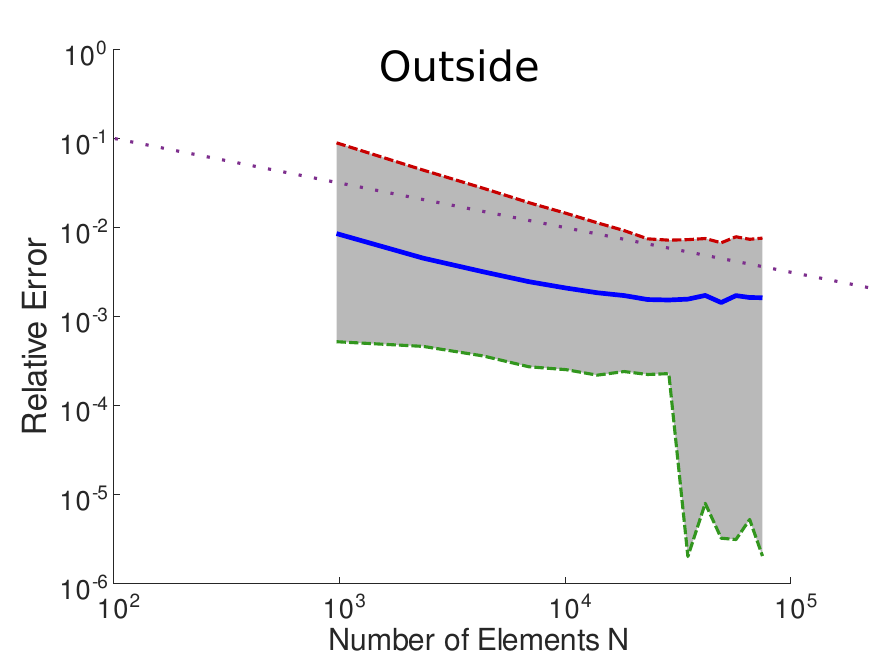}
    \end{tabular}
  \end{center}
  \caption{Relative errors as a function of the number of elements $N$ for the sphere benchmark at $f = 200$\,Hz. Left side: Error on the surface of the sphere. Right side: Error around the sphere. The continuous (blue) line depicts the mean error over all nodes, the gray area depicts the area between the maximum and the minimum  error (dashed lines). The dotted line depicts the function $\frac{1}{\sqrt{N}}$. }\label{Fig:Error200} 
\end{figure}
\begin{itemize}
\item For the integral over the boundary elements, the default maximum order of the Gauss quadrature is six. For regular integrals, this is a good compromise between accuracy and efficiency. However, sometimes this maximum order can be too small for quasi-singular integrals, even after applying the subdivision procedure (see Sec.~\ref{Sec:Quadrature}). %One option to deal with this problem, is to subdivide the element if the maximum Gauss order is reached. 
\item The default factor determining the multipole expansion length is set to 1.8 (see Sec.~\ref{Sec:Expansion}, Eq.~\ref{Equ:FMMTruncation}). In most cases this default value is sufficient for an accurate but fast calculation. However, if a higher accuracy is needed, this value needs to be increased at the price of a longer calculation time.
  % a factor of 12 is about single precision
\end{itemize}

\begin{figure}[!h]
  \begin{center}
    \begin{tabular}{rcrc}
      \raisebox{0.35\textwidth}{a)\hspace{-20pt}} &
      \includegraphics[width=0.46\textwidth]{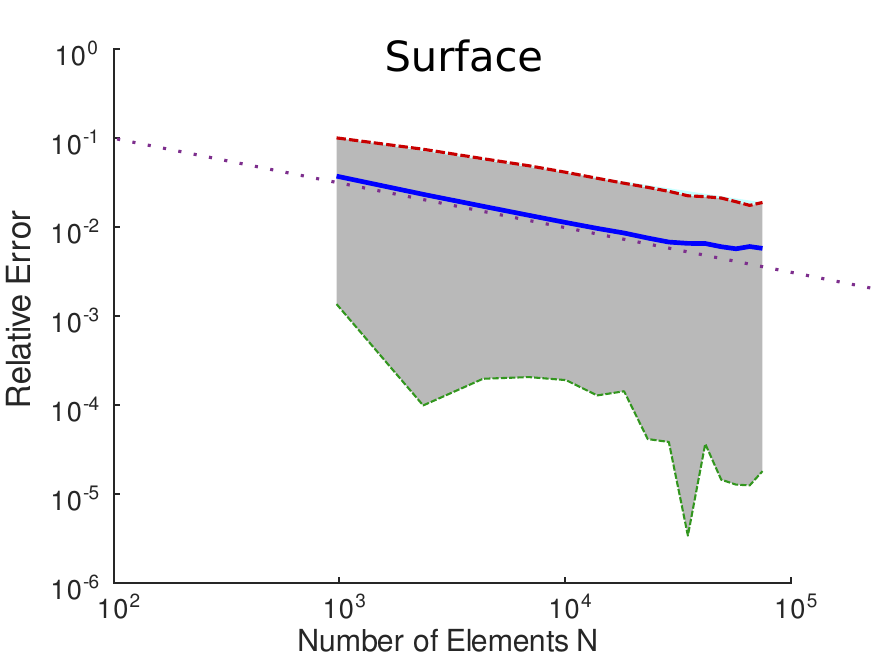} &
      \raisebox{0.35\textwidth}{b)\hspace{-20pt}} &                
      \includegraphics[width=0.46\textwidth]{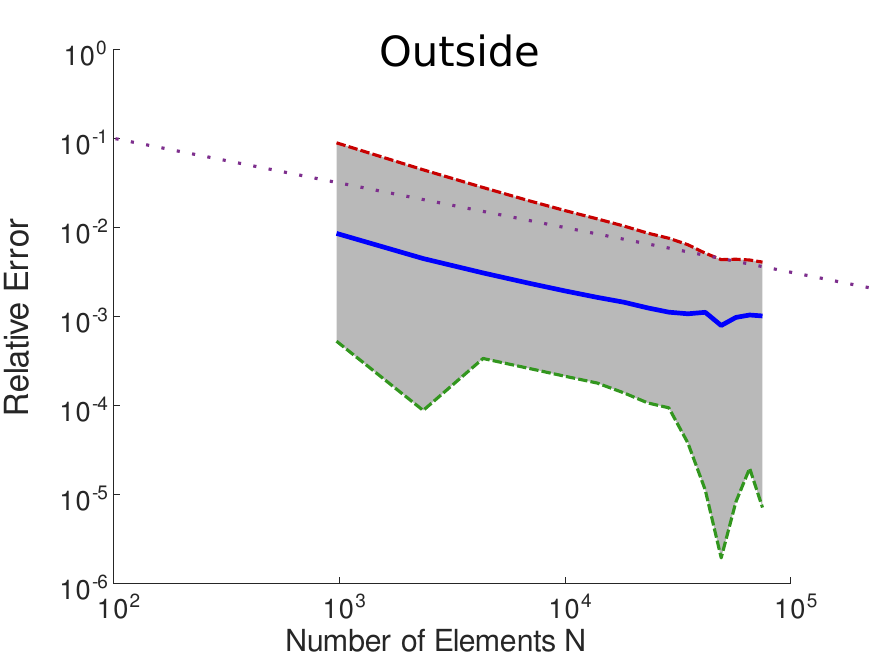}
    \end{tabular}
  \end{center}
  \caption{Relative errors with a higher accuracy setting as a function of the number of elements $N$ for the sphere benchmark at $f = 200$\,Hz. Left side: Error on the surface of the sphere. Right side: Error around the sphere. The continuous (blue) line depicts the mean error over all nodes, the gray area depicts the area between the maximum and the minimal error (dashed lines). The dotted line depicts the function $\frac{1}{\sqrt{N}}$.}\label{Fig:Error200rw15}
\end{figure}

In order to demonstrate the effect of the maximum order and the factor determining the multipole expansion length, we repeated the calculations, however, with two modifications of the parameters. First, the quasi-singular quadrature was modified such that each element was additionally subdivided if the estimated quadrature order was larger than 6. Second, the factor determining the multipole expansion length was set to 15 (see Eq.~\ref{Equ:FMMTruncation}). This resulted in slightly higher multipole lengths $L$. For example, for a mesh with $N = 74892$ elements, the truncation $L$ calculated for the multipole levels of $\ell = 1, 2$, and $3$  was $L_{1,2,3} = 12, 10, 9$ (in opposite to $L_{1,2,3} = 8, 8, 8$ obtained with the default setting), respectively. Fig.~\ref{Fig:Error200rw15} shows the relative errors obtained with the modified parameter settings. As expected, the errors decreased showing a better accuracy of the calculations. Note that the additional accuracy was paid with the longer computation time. For example, for a mesh with $N = 34992$ elements, the time to set up the linear system increased by factor \emph{two to three} as compared to that obtained with default parameters.

This observation has several other implications. First, it is possible to achieve a higher accuracy in the calculations. This is however linked with the price of longer calculations, even though the calculation time might be still reasonable. In our example, the calculation time was still within 5 minutes even for the mesh with $N = 74892$ elements. Further, users should ask themselves if the higher accuracy is actually required because the gain can be small (compare Fig.~\ref{Fig:Error200rw15} vs Fig.~\ref{Fig:Error200}). Optimizing for speed and accuracy is an issue and future work will include the option to automatically resample a fine mesh at low frequencies.
\subsubsection{Effect of the FMM}\label{Sec:NonZero}
In this section, we analyzed the effect of the FMM on the accuracy by calculating the relative difference
$$
\frac{|p_{\text{FMM}} - p_{\text{w/o FMM}}|}{|p_{\text{w/o FMM}}|}
$$
between the BEM calculations performed with and without the FMM. In the latter condition, we used a direct solver, i.e., the \texttt{zgetrf} routine from LAPACK~\cite{Andersonetal99} and limited $N$ to 8748 elements to keep the matrices at feasible sizes. The calculations were done for the frequencies 200 and 1000\,Hz.
\begin{figure}[!h]
  \centering
  \begin{tabular}{rcrc}
    \raisebox{0.35\textwidth}{a)\hspace{-20pt}} &
    \includegraphics[width=0.44\textwidth]{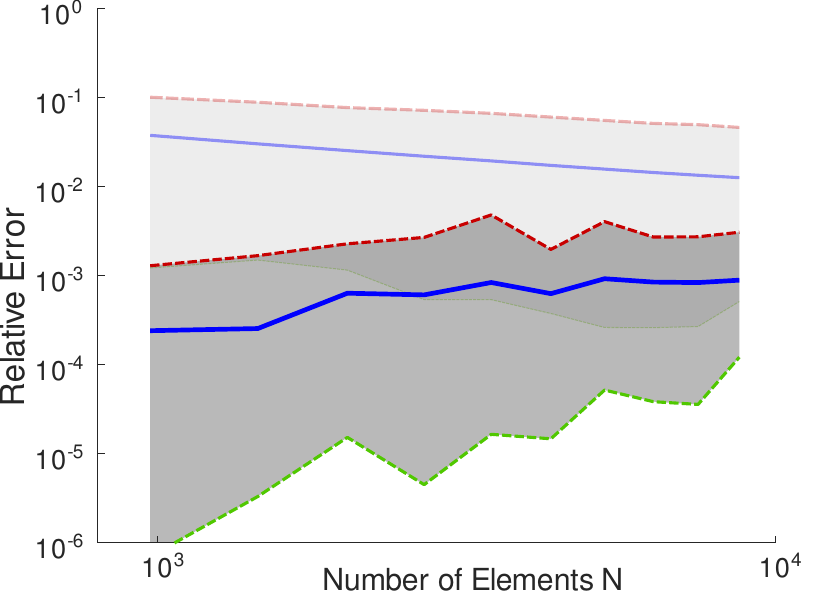} &
    \raisebox{0.35\textwidth}{b)\hspace{-20pt}} &
    \includegraphics[width=0.44\textwidth]{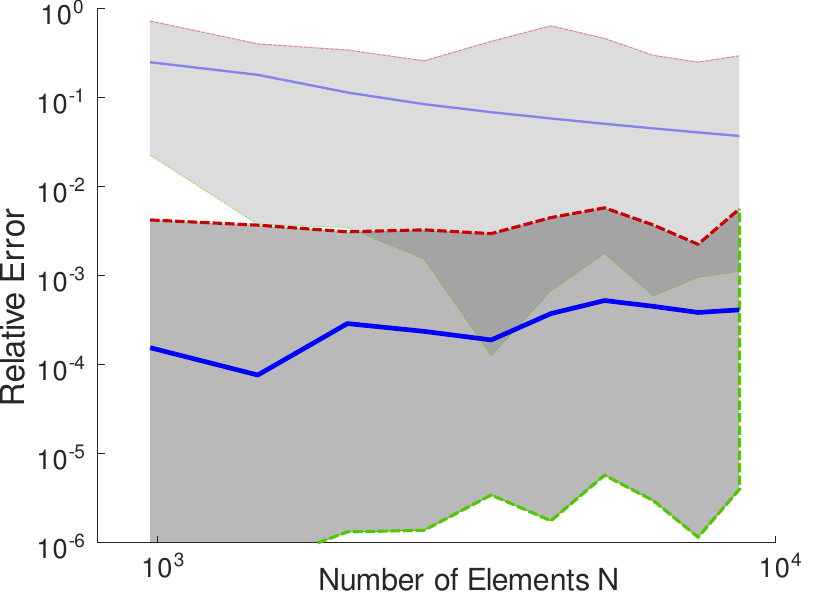}
  \end{tabular}
  \caption{Relative error between BEM and analytic solution (lighter colors), and the relative difference of the BEM calculations with and without FMM at 200\,Hz (a) and 1000\,Hz (b) as function of the number of elements $N$. The dashed lines depicts the maximum and minimal relative difference between conventional BEM and FMM solution, the continuous (blue) lines the mean distance between both solutions.}\label{Fig:FMMvsDirect}
\end{figure}
Fig.~\ref{Fig:FMMvsDirect} shows the relative errors obtained in the two cases as a function of $N$ for the two frequencies. This figure also shows the overall errors compared to the analytic solution in lighter colors. The relative errors introduced by the FMM were in the range of $10^{-3}$ showing only a minor effect on the overall error. Note that at 1000\,Hz (Fig.~\ref{Fig:FMMvsDirect}b), the wavelength is about $0.34$\,m, which means that the 6-to-8 elements per wavelength is only fulfilled for meshes with more than $6000$ elements. This explains the large overall errors in Fig.~\ref{Fig:FMMvsDirect}b.
\subsubsection{Memory requirement: Moderate mesh size}
In this section, we analyzed the memory requirements as an effect of the clustering. To this end, we compared the estimated and actual number of non-zero entries in the FMM matrices (see Sec.~\ref{Sec:FMM-Matrices}). We used the cube-based meshes and considered calculations for the frequency of $200$\,Hz for which we investigated the number of non-zero elements in the SLFMM and the MLFMM. The unit sphere was represented by a mesh with $N = 4332$ elements, for which the minimum, average, and maximum edge lengths was approximately 0.05\;m, 0.09\;m, and 0.15\;m, respectively. %These values correspond approximately to 21, 19 and 11 times the wavelength at  200\,Hz.
The 6-to-8 elements per wavelength rule is easily fulfilled in this example. 

\begin{figure}[!h]
  \begin{center}
    \includegraphics[width=0.45\textwidth]{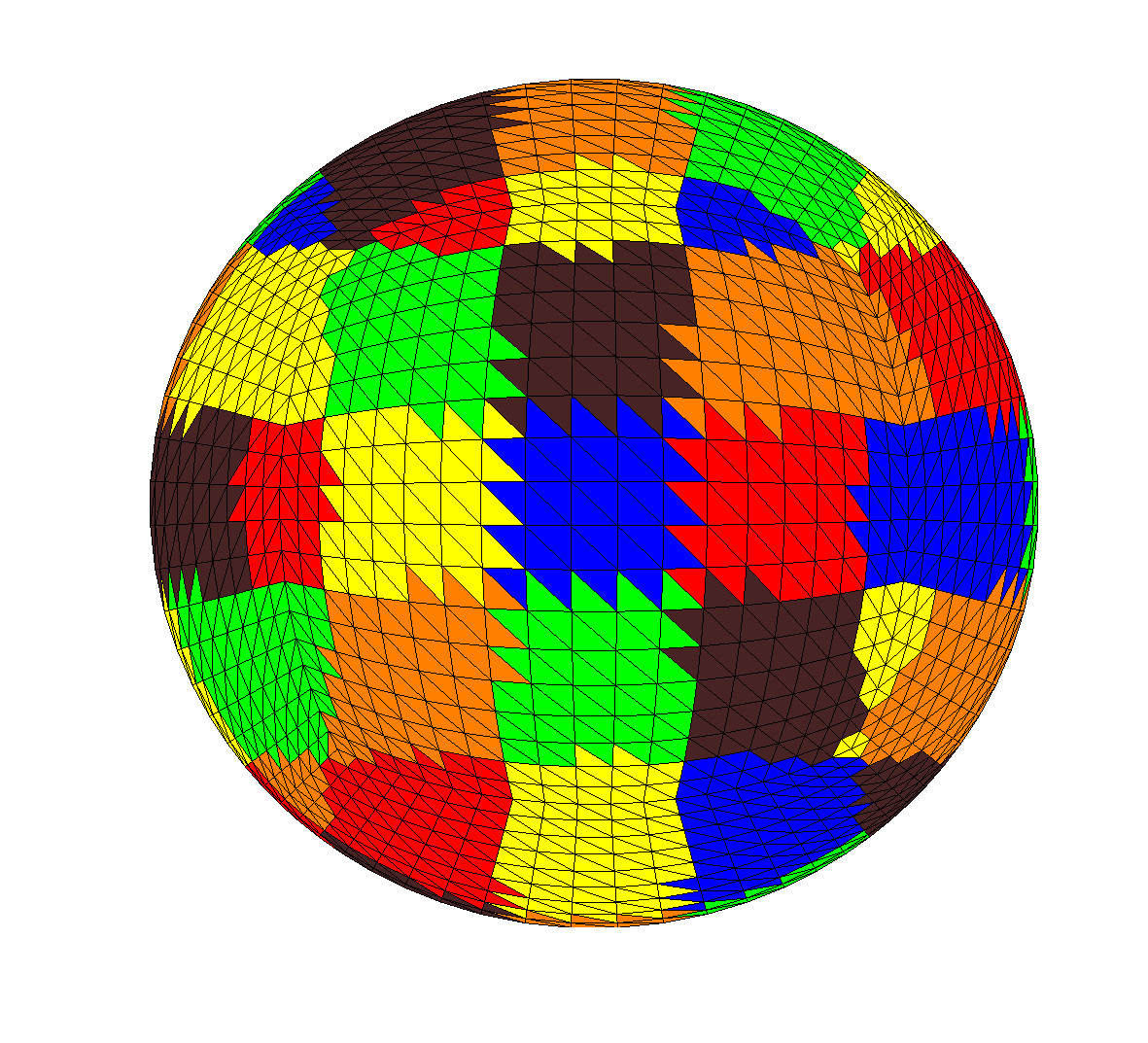}
    \caption{Example of default clustering for the unit sphere at 200\,Hz for a mesh with $N = 4332$ elements. The different colors/gray scales depict different clusters.}\label{Fig:Clustering1}
  \end{center}
\end{figure}

The clustering for the SLFMM resulted in $N_C = 90$ clusters with 29 and 72 elements. Fig.~\ref{Fig:Clustering1} shows the mesh and visualizes the obtained clusters. The average number of elements per cluster was $48$ elements, which is significantly different than the aimed estimation of $\sqrt N \approx 65$. This indicates that the actual clustering yields a different size even for a regular geometry such as the sphere. This difference has an impact on the estimation of the non-zero entries of the FMM matrices and therefor the accuracy of the estimation of the memory requirement.

Tab.~\ref{Tab:EntriesSLFMM_I} shows information about the mesh, the clustering for the SLFMM, the estimated (see also  \ref{App:Memory}), and the actual numbers of non-zero entries of the FMM matrices. It can be seen that the estimated number of non-zero entries of the midpoint-to-midpoint matrix $\bD$, which is dependent on a good estimation of the number of clusters, differs from the actual number of non-zero entries for this matrix. This is mainly because the estimation of the number of clusters was not correct. Nevertheless, if it is assumed that a complex valued entry has 16 Byte, this difference adds up to only about 7 Mbyte. Based only on the number of non-zero entries of the matrices, the overall memory used was roughly 75 Mbyte compared to about 300 Mbyte for the traditional BEM approach without FMM.
\begin{table}[!h]
  \begin{center}
    \begin{tabular}{ll}
      \multicolumn{2}{c}{Cluster Info}\\
      \hline
      Nr. Boundary Elems & $N = 4332$\\
      Nr. Clusters  & $N_C = 90$\\
      Truncation & $L = 8$\\
                   &\\
    \end{tabular}\hspace{1cm}
    $
    \begin{array}{crr}
      \text{Matrix} & \text{Actual} & \multicolumn{1}{c}{\text{Estimate}}\\
      \hline
      \bN & \num{2 700 492} &  \num{2 851 323}\\
      \bT,\bS & \num{554496} & \num{554496}\\
      \bD & \num{907520} &  \num{470249}\\
                 %   &        &  2L^2 N_C(N_C-10) = \num{921600}\\
    \end{array}
    $
 \end{center}
 \caption{Basic information about the mesh and the clustering for the SLFMM example and the estimated and actual number of non-zero entries of the FMM matrices.}\label{Tab:EntriesSLFMM_I}
\end{table}

The clustering for the MLFMM yielded $\ell_{\max{}} = 2$ levels. At the root level, we obtained $N^{[1]}_C = 56$ clusters, with the smallest, average, and maximum number of elements in an individual cluster being 30, 73, and 104 elements, respectively. The estimated number of clusters at the root level was $N^{[1]}_{C,\text{est}} = 0.9\sqrt{N} \approx 59$. At the leaf level, we obtained $N^{[2]}_C = 252$ clusters, with the minimum, average, and maximum number of elements in a cluster being 2, 17, 30 elements, respectively. The low number of elements demonstrates the effect of subdividing all parents on each level. On both levels, the expansion length was  $L_1 = L_2 = 8$. Tab.~\ref{Tab:EntriesMLFMM_I} shows the information about the clustering and the number of actual and estimated non-zero entries. For the estimation, it was assumed that each cluster has 10 clusters in its nearfield and that each parent has 4 children. The results show that in this case the estimates based on the number of elements $N$ are a good  approximation of the actual number of non-zero elements. The data presented in Tab.~\ref{Tab:EntriesMLFMM_I} suggests a memory consumption of the FMM matrices of about 75 MByte.

\begin{table}[!h]
  \begin{minipage}{0.4\textwidth}
    \begin{center}
      \begin{tabular}{ll}
        \multicolumn{2}{c}{Cluster Info}\\
        \hline
        Nr. Boundary Elems & $N = 4332$\\
        Nr. Clusters  & $N^{[1]}_C = 56$\\
                     & $N^{[2]}_C = 252$\\
        Truncation   & $L_1  = L_2 = 8$\\
      \end{tabular}
    \end{center}
  \end{minipage}
  \begin{minipage}{0.6\textwidth}
    \begin{center}
      $
      \begin{array}{crr}
        \text{Matrix} & \text{Actual} & \multicolumn{1}{c}{\text{Estimate}}\\
        \hline
        \bN & \num{931 200} & \num{1 083 000}\\
        \bT_\ell,\bS_\ell & \num{554496} & \num{554496}\\
        \bD_1 & \num{327 680} &  \num{370 048}\\
        \bD_2 & \num{1 193 984} &  \num{887 194}\\
      \end{array}
      $
     \end{center}
   \end{minipage}
   \caption{Basic information about the mesh and the clustering and estimated and actual number of non-zero entries of the MLFMM matrices.}\label{Tab:EntriesMLFMM_I} 
 \end{table}
\subsubsection{Memory requirement: Large mesh size}
In this section, we repeated the analysis from the previous example for a mesh with  $N = 8112$ elements. Tab.~\ref{Tab:EntriesSLFMM_II} shows the clustering information and the estimations of non-zero matrix elements for the SLFMM case. The clustering resulted in $N_C = 90$ clusters, with the number of elements in individual clusters ranging from $52$ to $144$ elements. Note that the estimation for the number of clusters $N_C = \sqrt{N}$ is quite accurate. When compared to the previous mesh with only half of the elements, it may be surprising that despite the higher number of elements the number of clusters did not increase. First, for the SLFMM the default clustering is determined by using $\sqrt{N}$ bounding boxes, which means that the difference in the number of bounding boxes for $N = 8112$ and $N = 4332$ is only about 25 boxes. Second, one has to distinguish between the bounding box and the cluster inside this box. It is possible that a bounding box contains no cluster or only small clusters that is merged with a neighboring cluster. Based on the data in Tab.~\ref{Tab:EntriesSLFMM_II}, the matrices for the SLFMM need about 186\,Mbyte compared to 1053 MByte for the BEM without FMM.
\begin{table}[!h]
  \begin{center}
    \begin{tabular}{ll}
      \multicolumn{2}{c}{Cluster Info}\\
      \hline
      Nr. Boundary Elems & $N = 8112$\\
      Nr. Clusters  & $N_C = 90$\\
      Truncation & $L = 8$\\
 %                  &\\
    \end{tabular}\hspace{1cm}
    $
    \begin{array}{crr}
      \text{Matrix} & \text{Actual} & \multicolumn{1}{c}{\text{Estimate}}\\
      \hline
      \bN & \num{8 678 432} & \num{7 306 206}\\
      \bT,\bS & \num{1 038 336} &  \num{1038336}\\
      \bD & \num{911 616} &  \num{923050}\\
    \end{array}
    $
  \end{center}
  \caption{Basic information about the clustering and the comparison between estimated and actual number of non-zero entries of the SLFMM matrices.}\label{Tab:EntriesSLFMM_II} 
\end{table}

Tab.~\ref{Tab:EntriesMLFMM_II} shows the clustering results and the estimations for the MLFMM. At the root level, the clustering yielded $\ell_{\max} = 2$ levels with $N^{[1]}_C = 90$ and $N^{[2]}_C = 360$ clusters. At the root level, the number of elements inside the individual clusters ranged from 52 to 144 elements and, on average, each cluster contained 90 elements. On the leaf level, the number of elements in the individual clusters were  between 4 and 54 elements with an average of 23 elements per cluster. At both levels, the expansion length of the MLFMM expansion was $L_1 = L_2 = 8$. The estimate for the nearfield matrix was slightly too small, but the overall estimates matched the actual numbers well. The memory needed by the FMM-matrices for this case is about 147 Mbyte.

\begin{table}[!h]
  \begin{minipage}{0.4\textwidth}
    \centering
    \begin{tabular}{ll}
      \multicolumn{2}{c}{Cluster Info}\\
      \hline
      Nr. Boundary Elems & $N = 8112$\\
      Nr. Clusters & $N^{[1]}_C = 90$\\
                   & $N^{[2]}_C = 360$\\
      Truncation & $L_1 = L_2 = 8$\\
    \end{tabular}
  \end{minipage}
  \begin{minipage}{0.6\textwidth}
    \centering
    $
    \begin{array}{crr}
      \text{Matrix} & \text{Actual} & \multicolumn{1}{c}{\text{Estimate}}\\
      \hline
      \bN & \num{2 403 648} & \num{2028000}\\
      \bT_\ell, \bS_\ell &  \num{1 038 336} & \num{1 038 336}\\
      \bD_1 & \num{911 616} & \num{736 128}\\
      \bD_2 & \num{1 747 968} & \num{1 661 337} \\
    \end{array}
    $
  \end{minipage}
  \caption{Basic information about the clustering and the comparison between estimated and actual number of non-zero entries of the MLFMM matrices.}\label{Tab:EntriesMLFMM_II} 
\end{table}
 
\subsubsection{Memory requirement: Bounding box length}
In this section, we investigated the effect of changing the initial length of the bounding box. The sphere was represented by $N = 8112$ elements as in the previous example, but the initial length of the bounding box at the root level was set to $0.7$\,m. For the estimation of the non-zero entries of the FMM matrices the actual number of clusters generated by \numcalc{} in the calculation were used. The calculations were done with the SLFMM for the cube-based meshes and the frequency of $200$\,Hz. 

Tab.~\ref{Tab:EntriesMLFMM_III} shows the clustering results and the estimations. The clustering generated $\ell_{\max} = 3$ cluster levels. The root level consisted of $N^{[1]}_C = 26$ clusters with 312 elements per cluster on average, ranging between 184 and 370 elements per cluster. On the second level $N^{[2]}_C = 144$ clusters were generated with 56 elements per cluster on average, ranging between 28 and 77 elements per cluster. On the leaf level $N^{[3]}_C = 536$ clusters were generated with 15 elements per cluster on average, ranging between 1 and 32 elements per cluster. On each level, the expansion length was $L_1 = L_2 = L_3 = 8$. The non-zero entries of the FMM matrices in this case would need about 173 MByte. 

This example illustrates that setting the bounding box length for the root level affects the clustering to some degree. It affected the number of levels, the number of elements per cluster, and the radii of each cluster. By adjusting the initial length of the bounding boxes at the root level, the radii of each cluster and in turn the length of the multipole expansion can be influenced. This helps in controlling the number of non-zero entries in the matrices $\bT_\ell$ and $\bS_\ell$ to some extent if $L$ becomes too large.

\begin{table}[!h]
  \begin{center}
  \begin{tabular}{lr}
    \multicolumn{2}{c}{Cluster Info}\\
    \hline
    Nr. Boundary Elems & $N = 8112$\\
    Nr. Clusters & $N^{[1]}_C = 26$\\
                      & $N^{[2]}_C = 144$\\
                      & $N^{[3]}_C = 537$\\
    Truncation & $L_1 = L_2 = L_3 = 8$
  \end{tabular}
  \hspace{1cm}
  $
  \begin{array}{crr}
    \text{Matrix} & \text{Actual} & \multicolumn{1}{c}{\text{Estimate}}\\
    \hline
    \bN & \num{1 525 972} & \num{1 225 411}\\
    \bT_\ell, \bS_\ell & \num{1 038 336} & \num{1 038 336}\\
    \bD_1 & \num{58624} &  \num{53248}\\
    \bD_2 & \num{641 024} &  \num{663552}\\
    \bD_3 & \num{2 334 864} &  \num{2 469 888}\\
   \end{array}$
     
   \end{center}
   \caption{Basic information about the clustering and the comparison between estimated and actual number of non-zero entries of the MLFMM matrices. The initial sub-box length was set to 0.7\;m, additionally it was assumed that each for each cluster the nearfield consists of 9 clusters.}\label{Tab:EntriesMLFMM_III}
\end{table}
\subsection{MLFMM: Head and Torso}\label{Sec:MLFMM4}
\begin{figure}[!h]
  \begin{center}
    \begin{tabular}{rcrc}
      \raisebox{0.35\textwidth}{a)} &                                     \includegraphics[height=0.35\textwidth]{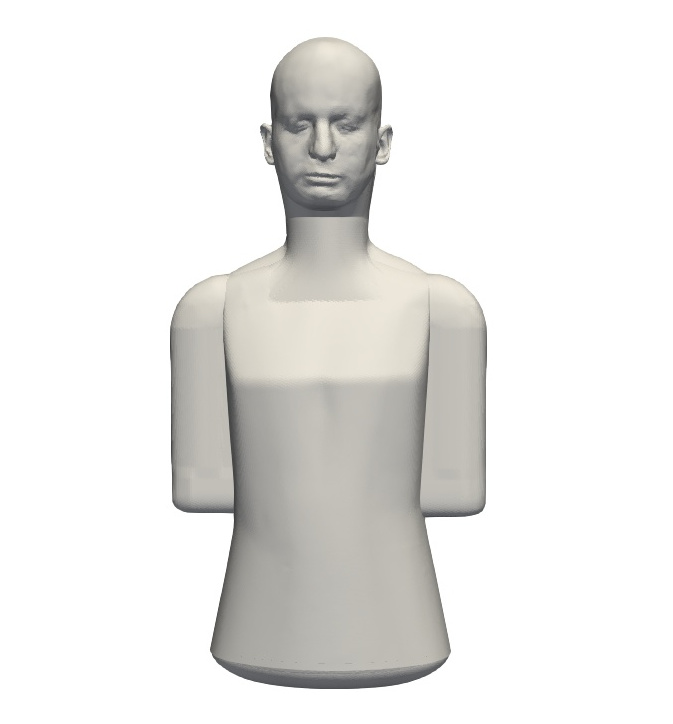} &                                       \raisebox{0.35\textwidth}{b)} &                                                         
 \includegraphics[height=0.35\textwidth]{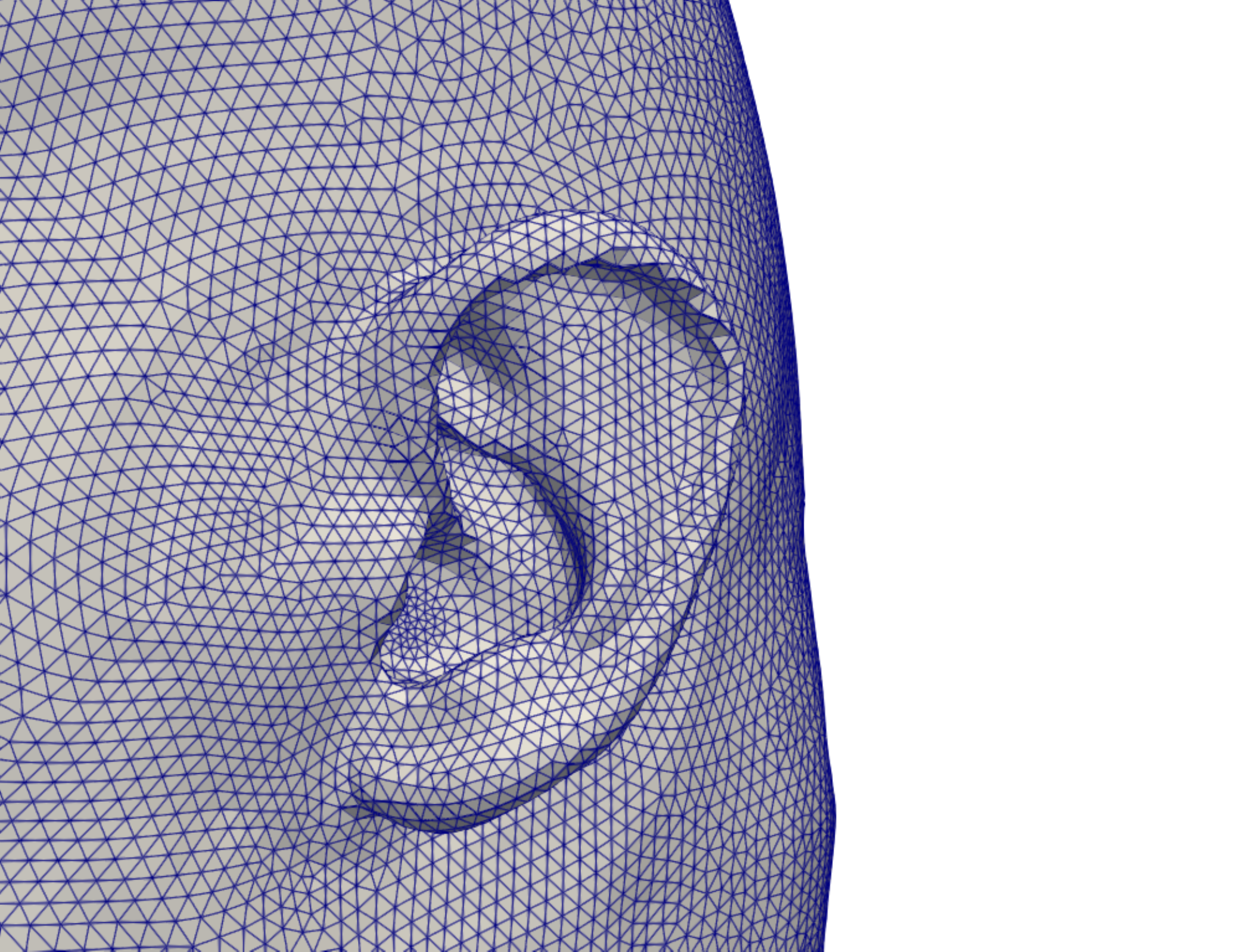}
    \end{tabular}
  \end{center}
  \caption{a) Geometry of the Head and Torso Example b) Mesh around the left pinna.}
\end{figure}

The example of a human head and torso of consisting of $\num{174 298}$ triangular elements was chosen to illustrate some of the restrictions on the current implementation of the FMM in \numcalc. The mesh is used to calculate the reflection and diffraction of sound waves at the human head, torso, and especially at the pinna. A point source was placed near the entrance of the pinna. There is also a large difference in the size of the elements between torso (relatively coarse discretization) and the human pinnae, where the discretization needs to be fine in order to correctly approximate the geometry and the solution. In this case, the minimum, average, and maximum edge length of all boundary elements were $0.7, 3.3,$ and $6.7$\,mm, respectively. The default clustering results in 3 levels with $N^{[1]}_C, N^{[2]}_C, N^{[3]}_C = 422, 1558, 6378$ clusters per level. The expansion length on each level was given by $L_1, L_2, L_3  = 31, 19$, and $11$, respectively. 

The amount of additional memory for storing the matrices $\bT_\ell$ and $\bS_\ell$ depends on $2L_\ell N$ on each level $\ell = 1,2,3$ and becomes quite large for this example (about 16\,GByte). The non-zero entries of $\bT_\ell$ and $\bS_\ell\; (\ell = 1,2)$ on levels $\ell = 1$ and $\ell = 2$ are about 61\% of all non-zero entries of the FMM matrices (cf. Tab.~\ref{Tab:EntriesMLFMM_V}). In this case, the gain in computing time by avoiding filtering/interpolation between levels does not outweigh the additional memory consumption. Of course, compared to the $O(N^2)$ memory needs of 486~GByte for the traditional BEM, the FMM is still very effective. The additional memory needed to explicitly store these matrices could be avoided, e.g, by filtering/interpolation algorithms for the FMM matrices. 

\begin{table}[!h]
  \begin{center}
    \begin{tabular}{lr}
      \multicolumn{2}{c}{Cluster Info}\\
      \hline
      Nr. Boundary Elems & $N = \num{170 914}$\\
      Nr. Clusters & $N^{[1]}_C = 403$\\
                & $N^{[2]}_C = 1512$\\
                & $N^{[3]}_C = 6226$\\
      Truncation & $L_1 = 37$\\
                & $L_2 = 20$\\
                & $L_3 = 16$
    \end{tabular}
    \hspace{1cm}
    $\begin{array}{cr}
      \text{Matrix} & \text{Entries}\\
      \hline
      \bN &  \num{74412824}\\
      \bT_1 = \bS_1 &  \num{335000756}\\
      \bT_2 = \bS_2 & \num{125843156}\\
      \bT_3 = \bS_3 & \num{42180116}\\
      \bD_1 & \num{333601540}\\
      \bD_2 & \num{36733916}\\
      \bD_3 & \num{51291416}
     \end{array}$
   \end{center}
   \caption{Basic information about the clustering and number of non-zero entries of the MLFMM matrices for the head-and-torso model.}\label{Tab:EntriesMLFMM_V}
\end{table}
\subsection{Duct problem}
\begin{figure}[!h]
  \begin{minipage}{0.44\textwidth}
    \centering
    \includegraphics[width=0.8\textwidth]{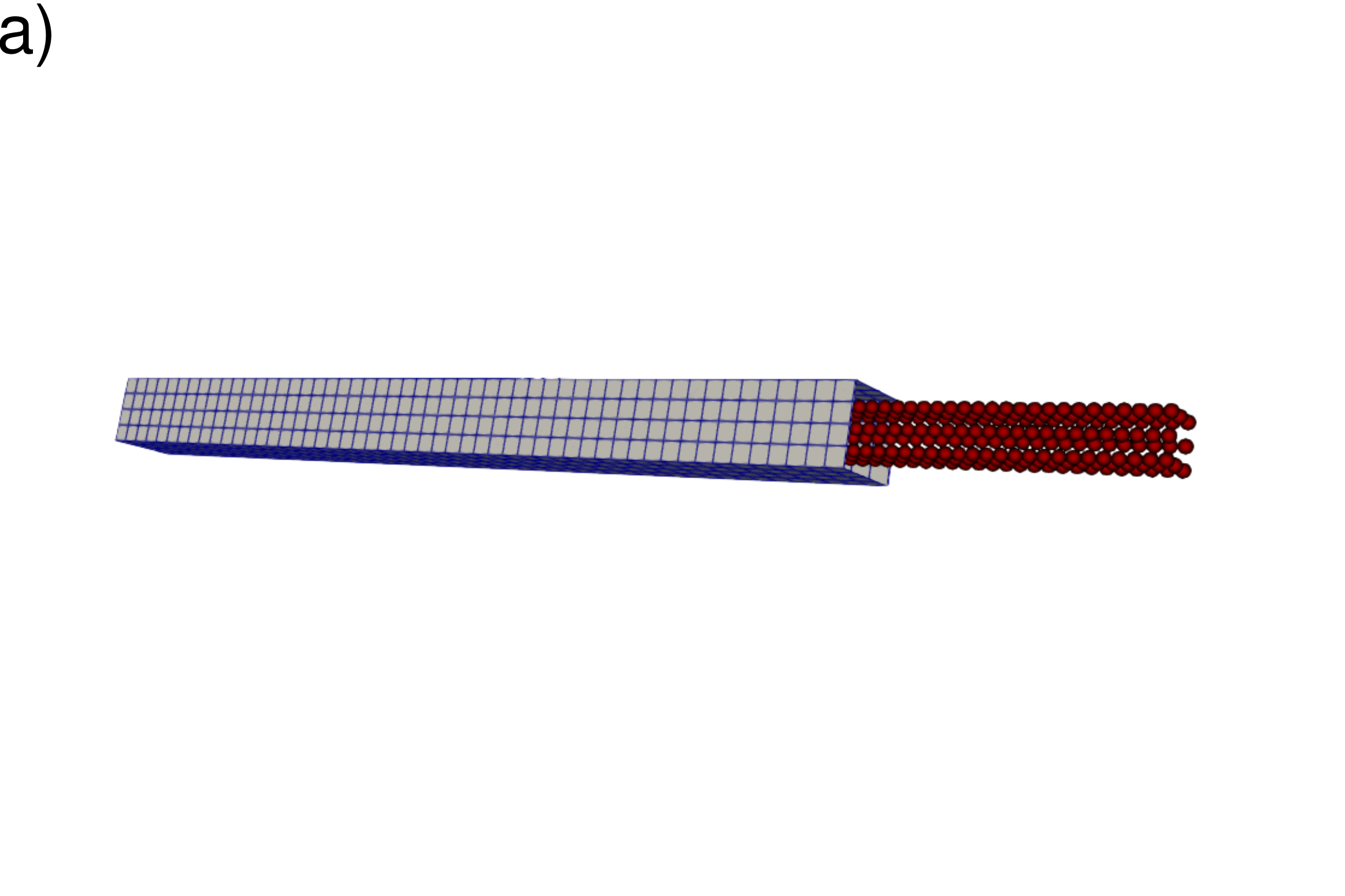}
  \end{minipage}
  \begin{minipage}{0.5\textwidth}
    \centering
    \includegraphics[width=0.98\textwidth]{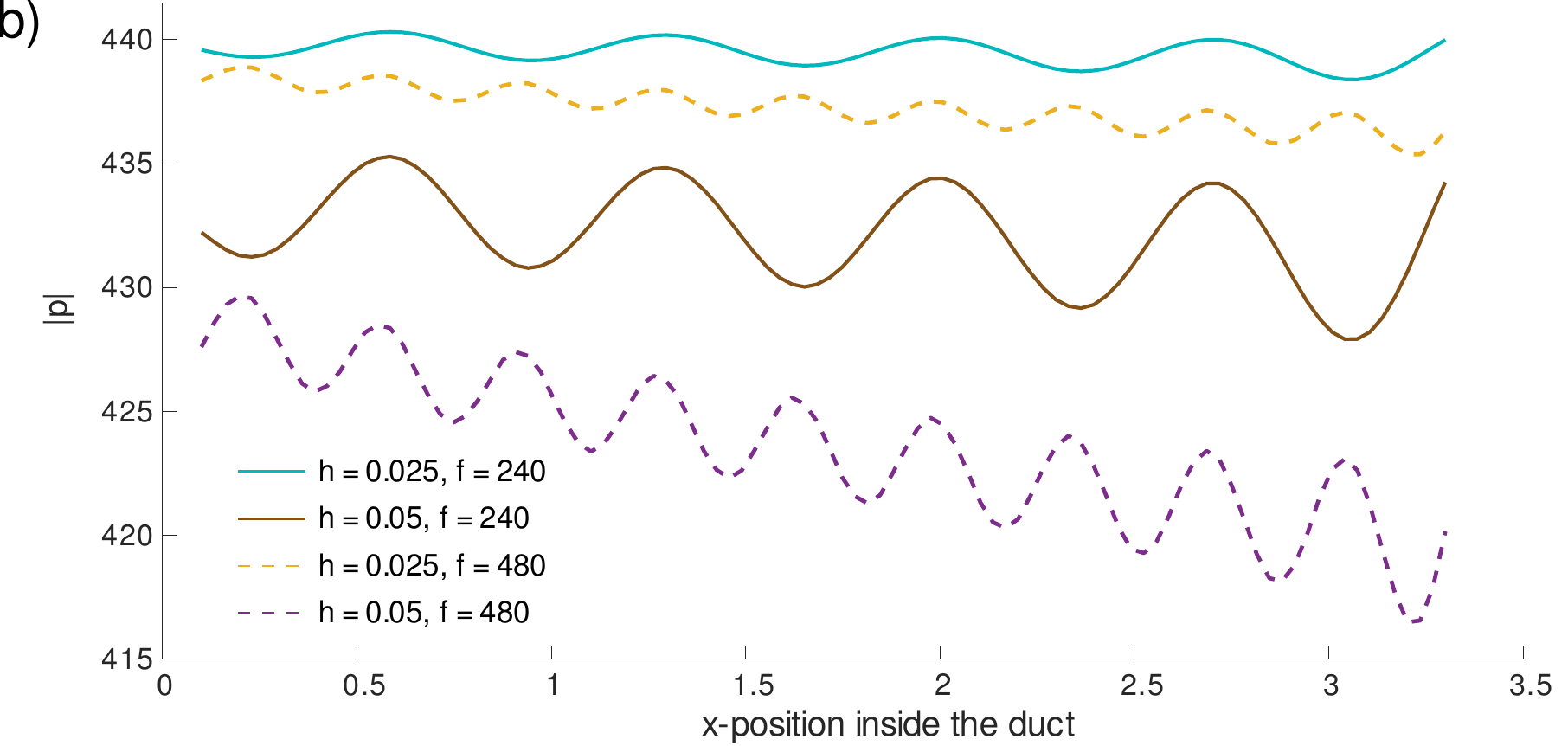}
  \end{minipage}
   \caption{a) Discretization of the duct and the evaluation points inside the duct. Only parts of the duct surface are displayed. b) Magnitude in Pa of the pressure along a line inside the duct for different element sizes $h$ at 240\,Hz and at 480\,Hz.}\label{Fig:InsideDuct}
\end{figure}

This problem is one of the benchmarks described in \cite{Hornikxetal15}. It describes the problem of calculating the acoustic field inside a closed sound-hard duct of length $L = 3.4$\,m that has a quadratic cross section with an edge length of $w = 0.2$\,m. At the closed side at $x = 0$, a velocity boundary condition with $v = 1$\,m/s is assumed, at the closed side at $x = L$, an admittance boundary condition with $Y = \frac{1}{\rho c}$ is assumed, where $\rho = 1.3$\,kg/m$^3$ and $c = 340$\,m/s are the density and the speed of sound in the medium, respectively. The quadratic cross section has the advantage that the geometry of the duct can be modelled accurately. A mesh for this problem including the positions of evaluation nodes inside the duct is partly depicted in Fig.~\ref{Fig:InsideDuct}a. The solution for these boundary conditions is given by a plane wave inside the duct traveling from one end to the other
\begin{equation}\label{Equ:PlaneWaveSol}
p(x) = - \rho c e^{\I kx},
\end{equation}
for the given speed of sound and the given density the magnitude of the solution Eq.~(\ref{Equ:PlaneWaveSol}) is $|\rho c| = 442$.

This problem has been investigated in several publications \cite{Kreuzer22,Marburg18} and it has been shown that calculations with collocation BEM result in a decaying field along the duct, where the rate of decay depends on the frequency and the element size along the duct.
As an example, the absolute value of the BEM solutions inside the duct for two different frequencies ($f_1 = 240$\,Hz and $f_2 = 480$\,Hz) and two different discretizations with edge lengths $h = 0.025\,$m and $h = 0.05$\,m are presented in Fig.~\ref{Fig:InsideDuct}b. At these frequencies, a smooth and robust solution is expected. It can be easily seen that the difference between the plane wave solution and the numerical solution grows along the length of the duct. Overall, the errors between calculated and analytic solutions are smaller than $6\%$. Nevertheless, this benchmark is one of the few problems where the 6-to-8-elements-per-wavelength rule may not be enough to provide results that are accurate enough.

\begin{figure}[!h]
  \begin{center}
    \includegraphics[width=0.46\textwidth]{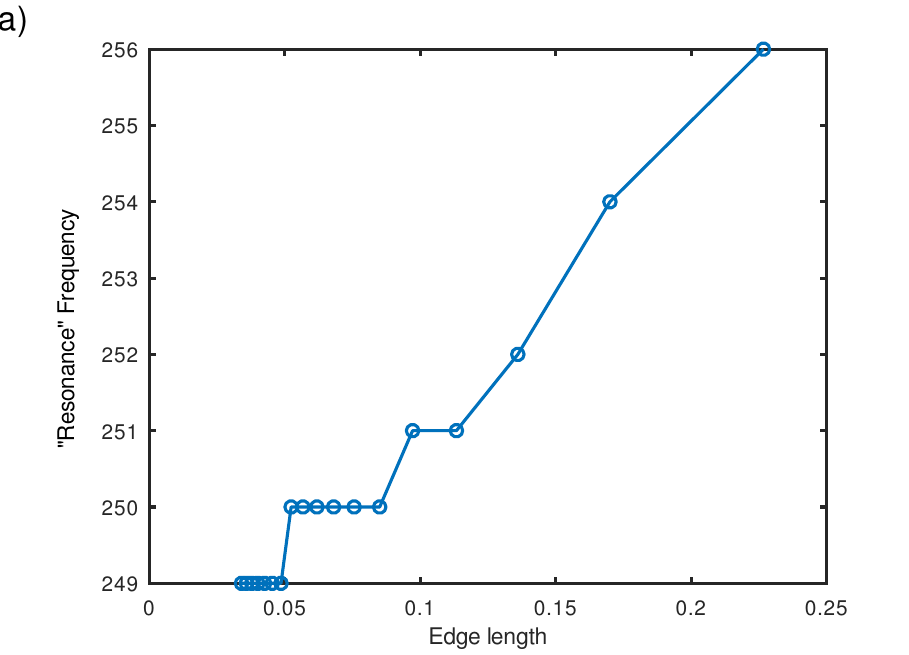}
    \includegraphics[width=0.46\textwidth]{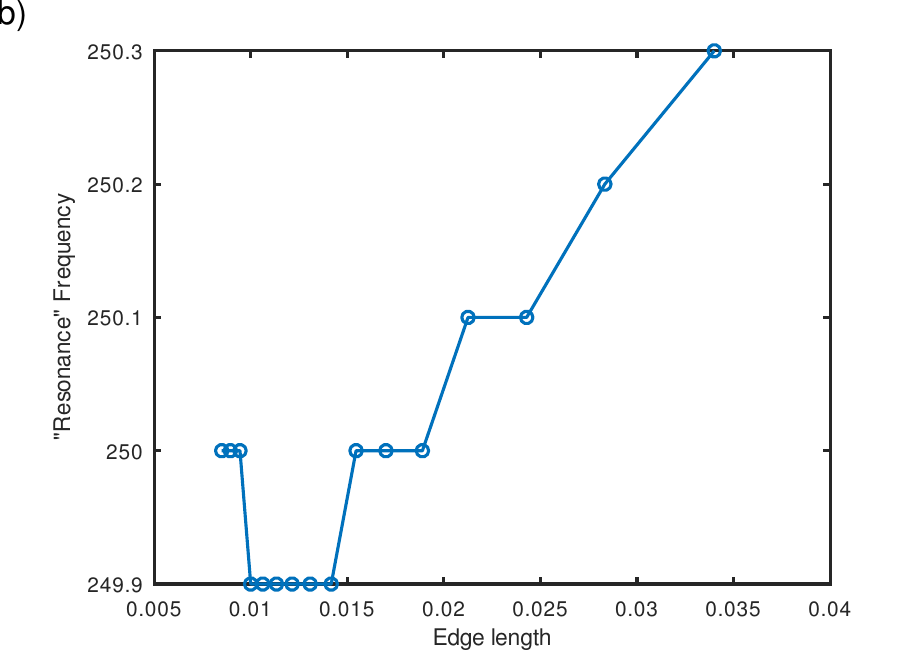}
    \caption{Frequencies at which the acoustic field inside the duct has maximum energy as a function of the element length $\ell$ in the $x$-direction. The whole duct is assumed to be sound hard. The figure on the right side is a zoomed version of the figure on the left.}\label{Fig:Resonances}
  \end{center}
\end{figure}

To better illustrate this fact, the boundary conditions at both ends at $x = 0$ and $x = L$ are changed to sound-hard conditions and an external source is placed inside the duct. For such a duct the theoretical resonance frequencies\footnote{When neglecting the diameter of the duct.} are at $f = 50\cdot n$\,Hz with $n \in \mathbb{N}$. Fig.~\ref{Fig:Resonances} shows the specific frequency around $250$\,Hz where the acoustic field inside the duct has maximum energy as a function of element size. It can be easily observed that the ``resonance'' frequency of the calculation is only close to the theoretical value if the edge length is about 0.1\,m, which is way shorter than suggested by the 6-to-8-elements-per-wavelength rule. At 250\,Hz, the wavelength is about $\lambda = \frac{340}{250}\approx 1.36$\,m, which would require an edge length of an element of approximately $h = 0.17$\,m if 8 elements per wavelengths are assumed. In Fig.~\ref{Fig:Resonances}, the frequency with maximum energy for this element size is about 4 to 5\,Hz away from the ``correct'' frequency.

On a side note: The frequency step size used for the computation in Fig.~\ref{Fig:Resonances} was 1\,Hz.  The calculated resonance around 250Hz is very sharp and due to numerical errors not at an integer frequency. The numerically calculated acoustic field  inside the duct is therefor bounded although it may become quite large. The values of 249\,Hz are a bit misleading because of the frequency step size. If a finer frequency grid is used around 250\,Hz (see Fig.~\ref{Fig:Resonances}b), we see that the calculated ``resonances'' are between 249.9\,Hz and 250\,Hz.

One explanation why 8 elements per wavelength are not enough to determine the ``correct'' resonance frequency may be that, in theory, the field inside the duct is a superposition of two plane waves travelling in opposite directions along the duct and that the zeros of the resulting standing wave at the resonance frequencies have to be resolved quite accurately.
\section{Conclusions}\label{Sec:Discussion}
In this article, the open-source BEM program \numcalc{}, which is part of the open-source software project \meshtohrtf, was described. The article aims at a balance between giving researchers detailed information about certain aspects of the BEM implementation (Sections~\ref{Sec:Implementation} and \ref{Sec:Benchmark}) and providing users with a general background of the BEM and the methods used in \numcalc{} (Sections~\ref{Sec:BriefBEM}, \ref{Sec:Critical}, and \ref{Sec:Benchmark}). The motivation behind the article was to provide users insights to decide whether \numcalc{} is a good tool for their problem at hand.

An important aspect of \numcalc{} to consider is that \emph{\numcalc{} does not check the mesh before calculations}. In cases of problematic meshes, the chance is a high that the program will terminate with an error message. In Sec.~\ref{Sec:Quadrature}, it was mentioned that if two elements are close together, one of them will be subdivided based on the distance of their midpoints. If elements are very irregular or even overlap, this distance criterion will never be fulfilled and the code terminates after a fixed number of subdivisions with an error message. However, in many cases, e.g., when the mesh contains small holes or normal vectors point to the wrong side, \numcalc{} will produce \emph{some} results, and \emph{it is up to the user to check the validity of the results.} Throughout the article general hints were added to help  users to see what can go wrong if non-ideal input data or parameters are used when running \numcalc{}. 
  
From our personal experience and bugreports/requests, we identified few types of common errors. The most common cause for such errors are badly shaped, loose or overlapping elements in the mesh.  A second possible problem is that the iterative solver does not converge. One reason for that can be distorted elements in the mesh or an expansion length $L$ of the FMM being too big. This can happen for meshes with a very large number of elements and if small and big elements are close to each other, in other words, if the elements in the mesh do not locally have the same size. Third, a large truncation parameter $L$ (see Section~\ref{Sec:Expansion}) may reduce the stability of the calculations. One potential solution may be a modification in the initial length of the boundary box in the FMM (see also Section \ref{Sec:FMM_I}).

A very common error is the wrong directions of the normal vectors. A normal vector in the wrong direction means that the role of interior problem and exterior problem is switched, which in turn means that in general sound sources are on the wrong side of the scatterer. This is one of the cases where \numcalc{} will finish without triggering an error, but the results will be completely wrong. Thus, the most important hint we can give is: \hint{Check the plausibility of your results.}

It is possible that problems arise due to a mesh reduction procedure: loose vertices may be inserted, elements may get twisted, normal vectors may be inverted, holes may appear in the surface or external sound sources may end up too close to the surface of the scatterer\footnote{Remember the Green's function has a singularity at 0, so does the Hankel function used in the FMM.}. Technically, most of these irregularities are not really problems for the code itself and no error warning will be triggered. Such errors can only be discovered by checking the results of the computation.
  
In Sec.~\ref{Sec:Benchmark}, we discussed the accuracy and the memory requirements of \numcalc. This section also provides users with a rough guideline on what to expect when using \numcalc{}. It reminds users on their responsibility to think about whether \numcalc{} is suitable for their problem at hand. In one of the examples we showed that in combination with the default FMM setting it is possible that the accuracy of the calculations may suffer at low frequencies when using fine grids (see Fig.~\ref{Fig:Error200}). This may be partly caused by a truncation parameter $L$ (see Eq.~\ref{Equ:FMMTruncation}) that was too low. We showed that this problems can be partly avoided by using a larger truncation parameter $L$ (see Fig.~\ref{Fig:Error200rw15}). However, in the current version of \numcalc{}, the parameter defining $L$ is hard coded, and the source code needs to be modified to change that. Before doing that, users should ask themselves, if it is really necessary to use such a fine grid for low frequencies in the first place. Naturally, a coarser mesh means less accuracy, but users should ask themselves if it is necessary to invest much computing time to achieve an accuracy of, e.g., $10 ^{-5}$ if other parameters like impedance or geometry have an uncertainty of about $10^{-2}$?

Section~\ref{Sec:Benchmark} also showed that the estimation of the memory requirements of \numcalc{} is not always straightforward as it depends on a good estimation of the size and number of the clusters at the different levels. To provide a better estimation for the required memory, \numcalc{} has an additional option to perform the clustering only (without the BEM computation, see \ref{Sec:MemoryEstim}). Based on the information about the size of the clusters, a rough estimation of the necessary RAM can be given.

The head-and-torso example in Sec.~\ref{Sec:MLFMM4} showed some of the limitations of \numcalc{} with respect to the required memory. As the length of the multipole expansion was high in the root level, the memory required to store the element-to-cluster matrices was also high. One way to avoid that is to restrict the size of the initial bounding box at the root level. This, in turn, results in clusters with smaller radii and smaller expansion lengths $L$. As the number of non-zero entries of $\bT$ and $\bS$ depends on $L$, memory consumption can be reduced this way to some degree.

Our discussion and examples provide a starting point for future improvements or adaptations to different computer architectures. If, for example, an architecture is used where memory is a bottleneck, e.g., GPUs, the memory consumption can be reduced by changing the clustering. However, as it can be seen from the results for the example with head and torso, the memory consumption for large problems demand more sophisticated approaches than the ones currently implemented in \numcalc. \numcalc{} was originally aimed at midsize problems where the number of elements are between $\num{60000}$ and $\num{120000}$. Future developments are planned to include dynamic clustering, interpolation and filtering between multipole levels to avoid the storage of the matrices $\bT$ and $\bS$ at levels different from the leaf level. Also, different quadrature methods for the integrals over the unit sphere  will be investigated.

With respect to accuracy, we also plan on using methods to adapt the mesh size to the frequency used, which would shorten computation time and enhances stability and accuracy in lower frequency bands. The simple example of plane wave scattering from a sound-hard sphere has already shown the general problem of having very fine grids at low frequencies. As the acoustic field  at low frequencies changes only slightly, the values on the collocation nodes are almost the same, which has a negative influence on the stability of the system.
\section*{Acknowledgements}
We want to express our deepest thanks to Z.-S.~Chen, who was the original programmer of the algorithms used in \numcalc{} and H. Ziegelwanger, who cleaned up the original code and created the original version of \meshtohrtf{}.
  
This work was supported by the European Union (EU) within the project SONICOM (grant number: 101017743, RIA action of Horizon 2020).
\appendix
\section{The Mesh}\label{Sec:Mesh}
A mesh is a description of the geometry of scattering surfaces using \emph{non-overlapping} triangles or plane quadrilaterals. The vertices of these elements are called BEM-nodes. Their coordinates are defined via node lists that contain a unique node-number and the coordinates of each BEM-node in 3D (see \ref{Sec:NodeList}). The elements are provided using an element lists containing a unique element number and the node-numbers of the vertices for each element. These lists also provide information if the element lies on the surface of the scatterer (BEM-element) or away from the surface (evaluation elements that are used if users are also interested in the field around or inside the object) (see \ref{Sec:NodeList}). The normal vector of each element always needs to point away from the object which is achieved by ordering the vertices of each element counterclockwise. The mesh should not contain any holes and each edge must be part of exactly two elements. The geometry described by the mesh should not contain parts that are very thin in relation to the element size, e.g., very thin plates or cracks, see also Sec. \ref{Sec:GreenSgl}. As a general rule of thumb there should be at least 6 to 8 elements per wavelength, see also Sec. \ref{Sec:Errororder}. The elements of the mesh should be as regular as possible. It improves the stability of the calculations, if elements are (at least locally) about the same size. Non-uniform sampling is possible as long as the change in element size is gradually, e.g., when using graded meshes towards singularities like edges or jumps in boundary conditions. A special kind of grading was developed in~\cite{Ziegelwangeretal16} that used coarser grids away from regions with much detail. In these regions the 6-to-8-element rule was violated. This grading was possible because for the specific problem of calculating HRTFs, which was the application at hand, the geometry of the pinna of interest is much more important than the rest of the head. Nevertheless, it is up to the user to have a good look at the result of the computation to judge its quality.

A specific boundary condition needs to be assigned to each element of the mesh: Pressure, particle velocity normal to the element (see \ref{Sec:BC}), or (frequency dependent) admittance, see also \ref{Sec:Admittance}. If no boundary condition is assigned to an element, it is automatically assumed to be sound hard, i.e., the particle velocity is set to 0 and the element is acoustically fully reflecting. Admittance boundary conditions are used to model sound absorbing properties of different materials. On each element, it is possible, to combine an admittance boundary condition with a velocity boundary condition. This additional velocity is used to model the effect of a vibrating part of the object that generates sound.
\section{Memory estimation}\label{App:Memory}
\subsection{Memory estimation SLFMM}\label{Sec:SLFMMEstim}
The number of non-zero entries of the nearfield matrix $\bN$ in Eq.~(\ref{Equ:SystemSLFMM}) is given by the product of the number of interactions between the elements of two clusters, the number of possible nearfield clusters for a cluster, and the number of all possible clusters. As every cluster has about $\sqrt{N}$ elements, a single nearfield cluster-to-cluster interaction adds  $N$ non-zero entries to the nearfield matrix. The number of clusters tn the SLFMM is about $\sqrt{N}$. In the worst case, each cluster has 27 clusters in its nearfield because a bounding box shares at least a face, an edge, or a vertex with 27 boxes (including itself). This implies that the nearfield matrix $\bN$ may have up to $27 N^{3/2}$ non-zero entries. However, one has to distinguish between the 3D bounding box and the cluster itself (see also Fig.~\ref{Fig:ClustervsBox}). In practice, one can assume that each cluster has on average about 9 to 10 nearfield clusters. It is fair to assume that for most calculations the matrix has about $9 N^{3/2}$ to $10 N^{3/2}$ non-zero entries. This assumption works well for small clusters where the surface of the scatterer inside the cluster has only little curvature. The cluster in Fig.~\ref{Fig:ClustervsBox}, for example, would probably have 9 clusters in its nearfield, one for each edge and each vertex of the cluster patch plus the patch itself. However, especially at lower cluster levels with bigger cluster radii, the assumption of low curvature is not always correct. In this case, the number of nearfield cluster can become larger than 9 or 10.

The local element-to-cluster expansion matrices $\bT$ and $\bS$ both have exactly $2L^2 N$ non-zero entries, where $2L^2$ is the number of quadrature points on the sphere that depends on the truncation parameter $L$. Since each element lies in exactly one cluster, there is only one local interaction between element and cluster midpoint resulting in $N$ local element-to-cluster expansions. 

For generating the cluster-to-cluster interactions matrix $\bD$, it is assumed that each cluster has about 9 clusters in the nearfield. The number of non-zero entries of the midpoint-to-midpoint matrix is a product of  number of clusters, number of far-field clusters, and the number of quadrature nodes on the unit sphere. On average $\sqrt N(\sqrt N - 9) 2L^2$ non-zero entries can be expected in total. The accuracy of this estimation depends on the uniformness  of the discretization and the smoothness of the scattering object.
\subsection{Memory Estimation MLFMM}\label{Sec:MLFMMEstim}
We assume that the scatterer has a regular geometry without many notches and crests,that each cluster has about 9 to 10 nearfield clusters, and that each parent cluster has only about 4 to 5 children.

In contrast to the SLFMM, the non-zero entries of the nearfield matrix $\bN$  now depend on the number of clusters on the highest (leaf-)level and the (average) number of elements in the leaf clusters. The number of non-zero entries on the leaf level can be estimated by the product of the number of clusters at this level, the number of element-to-element interactions between two clusters, and the number of nearfield clusters. Therefor, the nearfield matrix $\bN$ has about  $\frac{N}{25}\cdot 25^2 \cdot 9 = 225 N$ entries. The local element-to-cluster expansion matrices $\bT_\ell$ and $\bS_\ell$ depend on the number of terms used in the multipole expansion and have exactly $2 L^2_\ell N$ non-zero entries each.

The number of non-zero entries of the midpoint-to-midpoint interaction matrix at the root level can be derived similar to the SLFMM. At higher levels, the number of non-zero entries for $\bD_\ell$ is given by the product of the number of clusters at the level, the number of quadrature nodes on the sphere, and the number of clusters in the interaction list of each cluster. Assuming that each cluster has about 4 children and about 10 neighbouring clusters (without the cluster itself), the number of non-zero entries for $\bD_\ell  \approx 4^{\ell-1} 0.9 \sqrt{N} \cdot 2 L^2_\ell \cdot 40,\ell > 1$. In Section \ref{Sec:NonZero}, we compared these estimates with the actual number of non-zeroes in the different multipole matrices for a simple example. This showed that the estimates can only be rough estimates because the estimation for the number of clusters on each level is rather vague.
\section{Example Input File}\label{App:Input}
In the following, we present an example input file for the benchmark problem of a plane wave scattering from sound-hard sphere. The plane wave moves along the $z$ axis and the acoustic field will be calculated at 10 frequencies between $100$ an $1000$\,Hz. The different parts of the input file will be described in the following subsections.

Any line in the input file that starts with the hash tag '\#'  is ignored by \numcalc{} and can be used for comments. The unit for the coordinates is meters, the speed of sound is given in m/s, and the density is given in kg/m$^3$.
\begin{lstlisting}[basicstyle=\ttfamily,caption=Example input file for \numcalc{},label=Lst:InpFile,captionpos=b]
Mesh2HRTF 1.0.0
Mesh for plane wave reflection on a sound hard sphere
##
## lines starting with '#' are ignored
##
##
## Controlparameters I, reserved for later use
##
0 0 0 0 7 0
##
## Controlparameter II: Frequency Steps
## (dummy, nSteps, stepsize, i0)
# do 10 frequency steps with stepsize 1 and start at s1 = i0 + stepsize
1 10 1.0 0.0
##
## Frequency Curve 
#
# (Curvenr, Nr _of_Nodes)
0 2
# nodes of the piecewise linear curve
0.0 0
10.0 1000
## 1. Main Parameters I 
## (nBGrps, nElems, nNodes, min. ExpLength, nLevelsSym,
##  dummy, dummy, FMM_type, Solver)
## FMM_type = {0,1,4} = {No FMM, SLFMM, MLFMM}
## Solver = {0,4} = {CGS, direct Gauss Elimination}
2 28524 27042 8 0 2 1 4 0
##
## 2. Main Parameters II
# (nplanewaves, npointsrcs, ext. vs. int. probl,
#  default box length, default preconditioning)
1 0 0 0.0e0 0 
##
## 3. Main Parameters III, not used (yet, anymore)
0 0 0 0
##
## 4. Main Parameters IV
## dSoundSpeed, dDensity, harmonic factor
3.400000e+02 1.300000e+00 1.000000e+00
######################################
#
#  Define the geometry of the mesh
#
######################################
NODES
# Name of the file(s) containing the nodal data
TheBEMNodes.txt
TheEvalNodes.txt
ELEMENTS
# Name of the file(s) containing the nodal data
TheBEMElements.txt
TheEvalElems.txt
######################################
#
# define the boundary conditions
# First Elem, Last Elem, Realpart, Curve, Imagpart, Curve
#
######################################
BOUNDARY
# every element is sound hard, thus the default VELO = 0 
# condition can be used, no input necessary
RETURN
#
# external sound sources
# planewave
# (Sourcenr, direction(x,y,z), realpart, curve, imagpart, curve
# pointsource
# (Sourcenr, source position (x,y,z), realpart, curve, imagpart, curve
#
PLANE WAVES
0 0.000000e+00 0.000000e+00 1.000000e+00 1.000000e+00 -1 0.000000e+00 -1
#
# if some frequency dependent admittance, or external sources are used, 
# define the CURVES now
#
#CURVES
# no frequency dependent boundary conditions necessary
POST PROCESS
##
## BEM end of input ----
END
##
## this line seems to be important otherwise you end up in an
## endless loop, or at least add a return after END
\end{lstlisting}
In the current version v1.1.1 of \numcalc{}, the input files needs to be named \texttt{NC.inp}.
\subsection{Frequency-defining curves}\label{Sec:FreqCurv}
\numcalc{} numerically solves the Helmholtz equation in the frequency domain for a given set of frequencies defined by using frequency-defining curves. Frequency-defining curves are piece-wise linear functions that map a frequency step to a frequency or a frequency-dependent value of a boundary condition. The user has to provide the nodes of these curves. An example: If one would want to calculate the sound field for 10 frequencies between 100\,Hz and 1000\,Hz in steps of 100\,Hz, a simple way to describe the curve is given in Lst.~\ref{Lst:HRTFfreqs1}. The first non-comment line consists of a dummy parameter (\texttt{dummy = 1}) followed by the number of frequency steps to be calculated (\texttt{nSteps = 10}) followed by the distance $h_i$ between two frequency steps (\texttt{stepsize = 1.0}). The last entry in this line is given by the start index $i_0$ of the curve (\texttt{i0 = 0.0}).

The next line contains the unique number of the curve (in case of a frequency definition, the curve number always needs to be \texttt{Curvenr = 0}) and the number of nodes defining the piece-wise linear function mapping the frequency step $n$ to the respective frequency $f_n$. In Lst.~\ref{Lst:HRTFfreqs1}, this map is defined as the line between the two nodes
$(S_1, F_1)$ = ($0.0, 0.0$\,Hz) and $(S_2, F_2) = (10.0, 1000.0$\,Hz). The frequency for the $n$-th frequency step is then given as
\begin{equation}\label{Equ:Freqcurve}
f_n = F_1 + \frac{F_2 - F_1}{S_2 - S_1}((i_0 + h_i\cdot n - S_1) = 0\,\text{Hz} + \frac{1000 \,\text{Hz}}{10}(0.0 + 1.0 \cdot n) = 100\,\text{Hz}\cdot n,\, n = 1,\dots, 10.
\end{equation}
\begin{lstlisting}[basicstyle=\ttfamily,caption=Definition of 10 uniform frequencies between 100\,Hz and 1000\,Hz,label=Lst:HRTFfreqs1,captionpos=b]
  ## Controlparameter II: Frequency Steps
  ## (dummy, nSteps, stepsize, s0)
  # do 10 frequency steps with stepsize 1 and start at s_1 = s0 + stepsize
  1 10 1.0 0.0
  ##
  ## Frequency Curve 
  #
  # (Curvenr,  Nr_of_Nodes)
  0 2
  # nodes of the piece-wise linear curve
  # (step, frequency)
  0.0 0.0
  10.0 1000.0
\end{lstlisting}
Also note that for this example the 6-to-8-elements-per-wavelength rule is \emph{not} fulfilled for the two highest frequencies. Nevertheless, \numcalc{} will do the calculations, however, a warning will be displayed.

\begin{figure}[!h]
  \begin{minipage}{0.48\textwidth}
    \begin{lstlisting}[basicstyle=\ttfamily,caption=Definition of a frequency-defining curve for the third octave band example,label=Lst:FreqCurv,captionpos=b]
      #dummy n-steps stepsize 
      1     24  1.0 0.0
      # Curvenr. Nr._of_Nodes
      0      9
      #define the curve
      0.0  35.481
      3.0  44.668
      6.0  56.234
      9.0  70.795
      12.0 89.125
      15.0 112.202
      18.0 141.254
      21.0 177.828
      24.0 223.872
    \end{lstlisting}
  \end{minipage}
  % &
  \qquad
  \begin{minipage}[c]{0.48\textwidth}
    \includegraphics[width=0.92\textwidth]{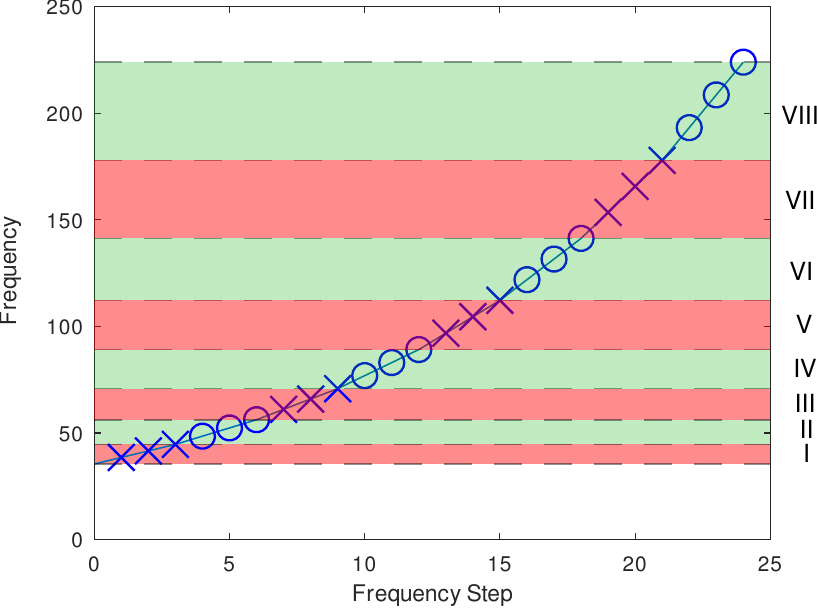}
    \caption{Frequency-defining curve for the third octave band.}\label{Fig:FreqCurvFig}
  \end{minipage}
\end{figure}

A more elaborate example is given in Lst.~\ref{Lst:FreqCurv}. This code snippet contains the definition of 8 third-octave frequency bands (see Fig.~\ref{Fig:FreqCurvFig}), where three frequencies are chosen in each band.\\
\subsection{Nodes and Elements}\label{Sec:NodeList}
The definition of the nodes (keyword \texttt{NODES}) and the elements (keyword \texttt{ELEMENTS}) are done by providing the names of the text files containing the data for the nodes and the elements. There needs to be at least one file for the nodes and at least one file for the elements, but it is also possible to use more than one text file to describe different parts of the mesh. This has the advantage that the surface of the scatterer and an evaluation grid can be defined independently as long as each node and each element have a unique identifying number. In the example input file Lst.~\ref{Lst:InpFile}, the information about the BEM nodes is given in \texttt{TheBEMNodes.txt}, whereas the coordinates of evaluation nodes not on the surface of the scatterer are given in \texttt{TheEvalNodes.txt}. The name of these files can be chosen arbitrarily by the user.

The structure of the node-defining files is as follows: The first line needs to contain the number of nodes in the file, the following lines contain a unique node number and the $x,y,z$ coordinates of the node in meters. The first few lines in \texttt{TheBEMNodes.txt} for this example looks like
\begin{lstlisting}[basicstyle=\ttfamily,caption=First few lines of \texttt{TheBEMNodes.txt}  containing the coordinates of the BEM nodes,label=Lst:TheBEMNodes,captionpos=b]
2168
1 -5.773503e-01 -5.773503e-01 -5.773503e-01
2 -5.975550e-01 -5.975550e-01 -5.346551e-01
3 -6.174169e-01 -6.174169e-01 -4.874348e-01
...
\end{lstlisting}
Jumps in the node numbers are allowed, which makes it easier to mix surface meshes and evaluation points.

The definition of the element files is as follows:
The first line contains the number of elements in the file. The next lines contain  a unique integer to identify the element and the node numbers of the vertices of the element followed by 3 integers. % denoting the type of the elements, zero and the group number of the element.
As an example, we show the first few lines of \texttt{TheBEMElements.txt} and \texttt{TheEvalElements.txt}
\begin{lstlisting}[basicstyle=\ttfamily,caption=First few lines of \texttt{TheBEMElements.txt} containing the definition of triangular surface elements,label=Lst:TheBEMElems,captionpos=b]
4332
1 381 382 458 0 0 0
2 381 458 457 0 0 0
3 382 383 459 0 0 0
...
\end{lstlisting}

\begin{lstlisting}[basicstyle=\ttfamily,caption=First few lines of \texttt{TheEvalElements.txt} containing the definition of the quadrilateral evaluation elements,label=Lst:TheEvalElems,captionpos=b]
24192
5000 2169 2170 2351 2350 2 0 1
5001 2170 2171 2352 2351 2 0 1
5002 2171 2172 2353 2352 2 0 1
...
\end{lstlisting}
The last three digits contain a flag for the type of elements, the number 0, and the group number of the element. \numcalc{} automatically determines from the number of the entries in each element line (7 or 8) if the element is triangular or quadrilateral (see Lst.~\ref{Lst:TheBEMElems} and \ref{Lst:TheEvalElems}). An element can either be a boundary element on the surface of an object (\texttt{type = 0}) or an element containing evaluation points (\texttt{type = 2}). Elements can also be assigned to different groups, but this feature is not recommended because the clustering for the FMM is done separately for each group. In general, it is sufficient to assign the BEM elements \texttt{GrpNr = 0} and evaluation elements \texttt{GrpNr = 1}.
\subsection{Boundary conditions}\label{Sec:BC}
Boundary conditions describe the acoustic behavior of the scattering objects. The definition section starts with the keyword \texttt{BOUNDARY}.  Three types of boundary conditions \{\texttt{PRES,VELO,ADMI}\} can be defined for each element. The general syntax for the boundary conditions is
\begin{verbatim}
ELEM e1 TO e2 bctype real(v0) curve1 imag(v0) curve2
\end{verbatim}
This line means that all elements with element numbers in the interval [\texttt{e1}, \texttt{e2}] have the boundary condition \texttt{bctype} $\in$ \{\texttt{PRES,VELO, ADMI}\} with value \texttt{v0}. \texttt{PRES} fixes the sound pressure for the element. If the pressure is set to zero, the element is sound-soft. \texttt{VELO} sets the particle velocity for the given element. If the velocity is set to zero, the element is sound hard. \texttt{ADMI} assigns an admittance to the element to model sound absorbing materials. In general, elements can only have one type of boundary condition, however, it is possible that elements have an \texttt{ADMI} as well as a \texttt{VELO} boundary condition. This is done to model (vibrating) sound sources on sound-absorbing surfaces. If no condition is assigned to an element, it is automatically assumed to be sound hard. 

As sound-absorbing behavior can be frequency dependent, \texttt{ADMI} conditions can be combined with a piece-wise linear curve (see Sec. \ref{Sec:FreqCurv}) that maps each frequency step to a scaling factor for the admittance. For the \texttt{ADMI} condition, \texttt{curve1} and \texttt{curve2} can contain the identifiers of piece-wise linear curves used to describe the frequency dependence of the real and imaginary part of the admittance. The curves are defined after the external sound source section. If no frequency dependence is needed, \texttt{curve1} and \texttt{curve2} are set to -1. The boundary condition section is closed using the keyword \texttt{RETURN}.
\subsection{Sound sources}\label{Sec:Sources}
Besides vibrating sound sources on the scatterer that are modelled by  \texttt{VELO} boundary conditions, sources away from the surface can be given by a combination of point sources and plane waves.

A point source $p_{\text{point}}(\bx) := P_0 \frac{e^{\I k ||\bx - \bx_0||}}{4\pi ||\bx - \bx_0||}$ at a source point $\bx_0$ is defined by its unique identifying  number, the coordinates of its origin \texttt{(X0,Y0,Z0)}, and its source strength \texttt{P0}:
\begin{verbatim}
POINT
Nr X0 Y0 Z0 Real(P0) curve1 Imag(P0) curve2
\end{verbatim}
A plane wave $p_{\text{planewave}}(\bx) := P_0 e^{\I k{\bf d}\cdot \bx}, ||{\bf d}|| = 1$ with direction \texttt{(DX,DY,DZ)} and strength \texttt{P0} is given by
\begin{verbatim}
PLANE
Nr DX DY DZ Real(P0) curve1 Imag(P0) curve2
\end{verbatim}
If the input is not dependent on the frequency, the curves identifiers are \texttt{curve1 = curve2 = -1}.
\subsection{Frequency definition curves for admittance}\label{Sec:Admittance}
If an admittance boundary condition or the strength of the external sound sources need to be frequency dependent, it is necessary to define additional scaling curves dependent on the frequency steps. The definition of the curves begins with the keyword \texttt{CURVES} followed by a line containing the \texttt{number\_of\_curves} and the largest number of nodes \texttt{nMaxCurves} that define these curves. Each single curve is defined by one line containing the unique curve number \texttt{icurve} and the number of nodes \texttt{inodes} defining the respective curve. The next \texttt{inodes} lines contain the nodes of the piece-wise linear curve similar to the frequency definition curve. An example:
\begin{verbatim}
CURVES
# number_of_curves max_number_of_points_per_curve
2 4
# curve_number  number_of_nodes_defining_curve
1 4
0.0 0.0 
1.0 1.0 
3.0 4.0
8.0 6.0
# curve_number number_of_nodes_defining_curve
2 3
0.0 0.0 
2.0 2.0
10.0 9.0
\end{verbatim}
defines two curves (one for the realpart, one for the imaginary part) that are used to define a frequency dependent admittance
\texttt{ELEM e0 TO e1 ADMI real(v0) 1 imag(v0) 2}, which assigns all elements between \texttt{e0} and \texttt{e1} a frequency dependent admittance
$
\alpha(s) = c_1(s)\, \text{real}(v0) + \I c_2(s)\, \text{imag}(v0) 
$
with 
$$
c_1(s) = \left\{
  \begin{array}{cl}
    s & \text{for } s\in [0,1],\\[2pt]
    1 + \frac32 (s - 1) & \text{for } s\in [1,3],\\[2pt]
    4 + \frac25 (s - 3) & \text{for } s\in [3,8],
  \end{array}
\right.
$$
and
$$
c_2(s) = \left\{
  \begin{array}{cl}
    s & \text{for } s \in [0,2],\\
    2 + \frac78(s-2) & \text{for } s\in [2,10].\\
  \end{array}
\right. 
$$
The frequency steps $s$ have already been defined via the frequency definition curve, i.e., the curve defined in Section \ref{Sec:FreqCurv}.
\section{Command-line parameters}\label{App:Commandline}
Besides setting parameters in the input file \texttt{NC.inp}, \numcalc{} offers additional command-line parameters.
\subsection{Estimation of memory consumption}\label{Sec:MemoryEstim}
The discussion about the non-zero entries of the multipole matrices in \ref{Sec:SLFMMEstim} and \ref{Sec:MLFMMEstim} has also shown that the clustering is very dependent on the mesh and that the numbers of clusters in the nearfield and the interaction list can vary to a large degree.
As it is hard to predict beforehand how the clustering will look like, \numcalc{} provides users with the option to just compute the clustering and the number of non-zero entries of the FMM matrices \emph{without} doing actual BEM calculations. If users start \numcalc{} with \texttt{NumCalc --estimate\_ram}, an estimation of the memory needed based on the entries of the multipole matrices is done. However, auxiliary variables and arrays to store mesh data and results are neglected in this estimation. Thus, the estimate slightly underestimates the actual memory needs. For the example of the sphere, \numcalc{} estimates the necessary  RAM to be 0.095 GByte, whereas Gnu/time \cite{Gnutime} listed the maximum memory consumption with 0.116 GByte. This difference can be explained partly by the number of evaluation nodes ($N_e = 24884$) used in this example. If $N_e = 70$, the memory consumption is reduced to about 0.1 GByte. This feature helps to get a rough idea to reserve the right amount of memory for calculations on a computer cluster.
\subsection{Parallelization}\label{Sec:Commandline}
As the BEM systems for different frequencies are independent of each other, parallelization is straight forward. \numcalc{} offers users the option to use the same input file for different ranges of frequency-steps. If, for example, only the first 10 frequencies should be calculated using one single thread, \numcalc{} can be started using \texttt{NumCalc -istart 1 -iend 10}, whereas a second thread can handle the next 10 frequencies using \texttt{NumCalc -istart 11 -iend 20}. However, it is up to the user to distribute the jobs over the different nodes of a computer cluster.
\section{Output}\label{App:Output}
The output of \numcalc{} can be split into two file types. First, there is \texttt{NC.out}\footnote{Or \texttt{NC\{step1\}-\{step2\}.out} if the command line parameters \texttt{-istart} and \texttt{-iend} are used.} containing general information, e.g., frequencies used, information about the mesh, the highest frequency where the 6-to-8-elements-per-wavelength rule is still valid and similar data.  This file also contains information about the FMM, e.g., the number of clusters, the expansion lengths for each level, information about the number of iterations, and the final error of the iterative solver.

Second, the results of the calculations are stored in the folder \texttt{be.out}. This folder contains several sub-folders \texttt{be.}$n$, one for each frequency step $n$. Inside each of these folders, one can find the files
\begin{itemize}
\item \texttt{pBoundary} containing the real and imaginary part of the sound pressure at the midpoints of each element,
\item \texttt{vBoundary} containing the real and imaginary part of the particle velocity normal to the surface at the midpoints of each element,
\item \texttt{pEvalGrid} containing the real and imaginary part of the sound pressure at each evaluation point,
\item \texttt{vEvalGrid} containing the real and imaginary part of the particle velocity in the $x$, $y$, and $z$ direction.
\end{itemize}
The structure of having separate folders for each frequency step has the advantage that the calculation for different frequency steps can be easier distributed on different nodes of a computer cluster.
%%%%%%%%%%%%%%%%%%%%%%%%%%%%%%
\bibliography{mesh2hrtf}{}
\bibliographystyle{elsarticle-num}
\end{document}